\documentclass[11pt]{amsart}
\pdfoutput=1
\usepackage{amssymb}
\usepackage{amsmath}
\usepackage{amsfonts}

\usepackage{graphicx}
\makeatletter 
 
 \@addtoreset{equation}{section}
\makeatother

\newtheorem{theorem}{Theorem}[section]
\newtheorem{lemma}[theorem]{Lemma}

\newtheorem{prop}[theorem]{Proposition}
\newtheorem{cor}[theorem]{Corollary}

\newtheorem{question}[theorem]{Question}

\theoremstyle{definition}

\newtheorem{remark}[theorem]{Remark}
\newtheorem{example}[theorem]{Example}

\newcommand{\on}{\operatorname}

\newcommand{\Spinc}{\on{Spin}^c}

\newcommand{\h}{\widehat}

\newcommand{\f}{\frac}
\newcommand{\ZZ}{\mathbb{Z}}

\newcommand\goth[1]{\mathfrak{#1}}
\newcommand{\s}{\goth{s}}

\newcommand{\sign}{\on{sign}}

\newcommand{\se}{\scriptstyle}

\begin{document}

\author{Shelly Harvey}
\address{Department of Mathematics, Rice University, Houston, TX 77005}
\email{shelly@rice.edu}
\author{Keiko Kawamuro}
\address{Department of Mathematics, Rice University, Houston, TX 77005}
\email{kk6@rice.edu}
\author{Olga Plamenevskaya}
\address{Deparment of Mathematics, Stony Brook University, Stony Brook, NY 11790}
\email{olga@math.sunysb.edu}
\title{On transverse knots and branched covers}
%\subjclass{?}

\begin{abstract} We study contact manifolds that arise as cyclic branched covers of transverse knots 
in the standard contact 3-sphere. We discuss properties of these contact manifolds and 
describe them in terms of open books and contact surgeries.
In many cases we show that such branched covers are contactomorphic for smoothly isotopic 
transverse knots with the same self-linking number. These pairs of knots include most of   
the non-transversely simple knots of Birman--Menasco and Ng--Ozsv\'ath--Thurston.

\end{abstract}
\maketitle

\section{Introduction}

In this paper, we consider transverse knots in $(S^3, \xi_{std})$, i.e. knots that are transverse to the 
contact planes of the standard contact structure $\xi_{std} =\ker(dz-ydx)$. 

A simple ``classical'' invariant is given by the self-linking number $sl$ of a transverse knot. 
However, if $L_1$, $L_2$ are two transverse knots  that are 
smoothly isotopic and share the same self-linking number,  $L_1$ 
and $L_2$ do not have to be transversally isotopic: this phenomenon 
was first discovered in \cite{EH1} and \cite{BM}, and more examples were recently obtained 
in \cite{NgOT}. 

The goal of this paper is to study  contact manifolds that arise as cyclic  covers
branched over transverse knots and links in $(S^3, \xi_{std})$. (Such cyclic covers carry natural contact structures lifting $\xi_{std}$.) The main question 
we would like to address is 

\begin{question}  Suppose that transverse knots $L_1, L_2$ are smoothly isotopic, and $sl(L_1)=sl(L_2)$.
Fix $p\geq 2$. Are $p$-fold cyclic covers branched over $L_1$ and $L_2$ contactomorphic?  \label{quest}  
\end{question} 
 
Finding two such non-contactomorphic covers would imply that the induced contact structure on the branched 
cover is an efficient invariant of transverse knots. On the other hand, a positive answer to the above question 
for any pair of knots means that the cyclic branched covers are insensitive to the subtler structure of transverse knots. 

While we found no examples of non-contactomorphic branched covers, we are able to answer Question \ref{quest} 
positively in many special cases. In particular, we show that
branched cyclic covers of any degree are contactomorphic for all examples of Birman--Menasco \cite{BM,BM1}, 
and that branched double covers for many examples of \cite{NgOT} are also contactomorphic.  
We prove 

\begin{theorem}  The $p$-fold cyclic branched covers  of transverse links $L_1$ and $L_2$ are contactomorphic for all $p$ if:

\begin{itemize} 
\item  $L_1=L^+$ and $L_2=\bar L^-$ are a positive and a negative transverse push-offs of a Legendrian link $L$ and its Legendrian mirror $\bar L$, 
 or 
\item $L_1$ and $L_2$ are given by transverse 3-braids related by a negative flype.  
\end{itemize} 
Moreover, the branched double covers are contactomorphic for arbitrary transverse braids related by a negative flype.
\label{introthm}
\end{theorem} 

In fact, we are able to prove a little more (see Section~\ref{distinguish}). We also note
that all examples of Birman--Menasco satisfy the second condition of Theorem~\ref{introthm}. 
(A {\em negative flype} is a braid move introduced in \cite{BM}; see also Figure~\ref{neg-flype} in Section~\ref{distinguish}.)

Let $\xi_p(L)$ denote the natural contact structure on the branched $p$-fold cyclic cover $\Sigma_p(L)$
of a transverse link $L$ as explained in Subsection~\ref{subsec-induced-str}.  We describe the contact manifolds $\Sigma_p(L)$ in two ways. First, we give an open book decomposition supporting $\xi_p(L)$. If $L$  
is represented as a transverse $n$-braid, an open book decomposition of $(\Sigma_p(L), \xi_p(L))$ can be obtained as a lift of the open book for $S^3$ whose binding is the braid axis, and the page is a disk meeting 
$L$ transversely at $n$ points.  The monodromy for the resulting open book can be read off the braid word. More precisely, a positive crossing in the braid  contributes $(p-1)$ positive Dehn twists to the monodromy, while a negative crossing contributes $(p-1)$ negative Dehn twists. Second, we give contact surgery diagrams \cite{DG,DGS} for these contact manifolds. We find that a positive (resp. negative) crossing in the braid corresponds
to a Legendrian surgery (resp. $(+1)$ contact surgery) on $(p-1)$ standard Legendrian unknots.
Interestingly, it turns out that the linking between these $(p-1)$ unknots depends on 
the sign on the crossing: while for a negative crossing the surgery is performed on 
{\em unlinked} unknots, for a positive crossing the unknots are linked. (This phenomenon arises in the smooth 
setting as well, where the construction can be thought of as a version of the Montesinos trick for higher order covers.) We refer the reader to Lemma \ref{Dehn twist lemma} and Theorem \ref{SurgThm} for precise statements.

This description yields a few properties of $p$-fold cyclic branched covers; we can determine whether they are  tight or overtwisted in certain cases  and describe the homotopy invariants of the contact structures.

\begin{theorem} \label{tight-ow} 
 
The contact manifold $(\Sigma_p(L),\xi_p(L))$ is Stein fillable if the transverse link
$L$ is represented by a quasipositive braid; it is overtwisted if $L$ is obtained as a transverse stabilization of another transverse link. 
\end{theorem}

In fact, in Section~\ref{properties} we  prove overtwistedness for a much wider class of contact structures.

\begin{theorem} \label{homotp} 
Fix $p\geq 2$. Let $\s_L$ be the $\Spinc$ structure induced by
$\xi = \xi_p(L)$. Then  $c_1(\s_L)=0$. The three-dimensional invariant $d_3(\xi)$ is
completely determined by the topological link type of $L$ and its
self-linking number $sl(L)$.
\end{theorem}

The present paper continues the work started by the third author in \cite{Pla}, 
where Question \ref{quest} was studied for the case of branched double covers. 
(The paper \cite{Pla} was written before the advent of Heegaard Floer transverse
invariants \cite{OST}, and the only explicit examples of non-transversely simple knots
available then were the 3-braids of \cite{BM}.) The techniques from \cite{Pla} are useful for 
the higher order covers as well; in particular, Theorems \ref{homotp} and \ref{tight-ow} are 
direct extensions of results of \cite{Pla}.

%%%%%%%%%%%%%%%%%%%%%%%%%%%%%%%%

\section{Preliminaries}

In this section, we fix notation and collect the necessary facts about transverse knots, open books 
and contact surgeries, referring the reader to \cite{Et, Et3, DGS} for more details. We assume  
that all 3-manifolds are closed, connected and oriented, 
and all contact structures are co-oriented.

\subsection{Transverse knots as braids} 

It will be helpful to represent transverse links by closed braids. For this, consider the symmetric version of the standard contact structure $(S^3,\xi_{sym})$ with $\xi_{sym}=\ker(dz+xdy-ydx)$. Then, any closed braid about $z$-axis can be made transverse to the contact planes; moreover, any transverse link is transversely isotopic to a closed braid in $(S^3,\xi_{sym})$ \cite{Be}. Equivalently, we can  consider transverse braids in the contact structure $\xi_{std} =\ker(dz-ydx)$, for example assuming that our braids are satellites of a fixed standard Legendrian unknot with $tb=-1$.

To define the self-linking number $sl(L)$, trivialize the plane field $\xi$,
and let the link $L'$ be the push-off of $L$ in the direction of the first coordinate
vector for $\xi$. Then, $sl(L)$ is the linking number between $L$ and $L'$.
Given a closed braid representation of $L$, we have
\begin{equation}\label{sl-braid}
sl(L)=n_+ -n_- - b,
\end{equation}
where $n_+$ ($n_-$) is the number of positive (negative) crossings, and $b$ is the braid index.

The {\em stabilization} of a transverse link represented as a braid is
 equivalent to the negative braid stabilization, i.e. adding an extra strand
and a negative kink to the braid. If $L_{stab}$ is the result of stabilization of $L$, then
\begin{equation}
sl(L_{stab})=sl(L)-2.
\end{equation}
Note that the positive braid stabilization does not change the transverse type of the link.

Abusing notation, we will often identify a transverse link with its braid word, 
writing it in terms of the standard generators $\sigma_1, \sigma_2,\cdots$ and their inverses.

Another useful way to think about transverse knots is as push-offs 
of Legendrian knots. Indeed, a given Legendrian knot yields two transverse knots 
(a positive and negative push-off), whose self-linking number is $tb(L)\pm r(L)$. 
This description is used in \cite{OST, NgOT}.

\subsection{Open books}

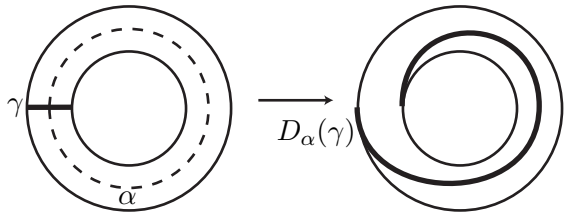
\begin{figure}[htpb!]
\begin{center}
\begin{picture}(203, 78)
%\put(0,0){\includegraphics[scale=0.4]{rh-twist}}
\put(0,0){\includegraphics{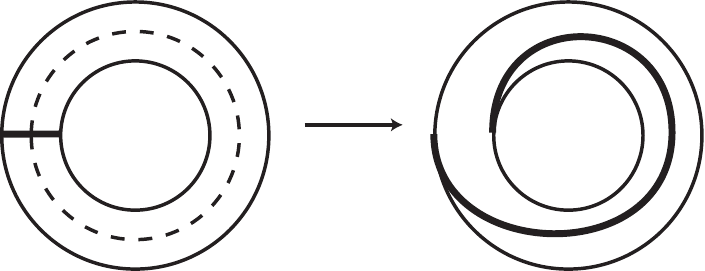}}
\put(35, 3){$\alpha$}  
\put(-7, 38){$\gamma$}
 \put(95, 27){$D_\alpha(\gamma)$}
\end{picture}
\caption{A right-handed Dehn twist $D_\alpha$ about $\alpha$. }\label{rh-twist}
\end{center}
\end{figure}

An open book decomposition of a 3-manifold $M$ is a pair $(S, \phi)$ of a surface $S$ with non-empty boundary $\partial S$ and a diffeomorphism $\phi$ of $S$ with $\phi|_{\partial S}=id$, such that 
$M \setminus \partial S$ is the mapping torus $S \times [0,1] / \sim $, where   $(x,1) \sim (\phi(x), 0)$. 
The surface $S$ is called a {\em page} and $\partial S$ the {\em binding} of the open book. 
By the celebrated work of Giroux \cite{Gi}, contact structures on $M$ up to an isotopy are in one-to-one 
correspondence with open book decompositions of $M$ up to stabilization. Stabilization of an open book 
consists of plumbing a right-handed Hopf band, i.e attaching a 1-handle to a page and composing the monodromy 
with a right-handed Dehn twist along an arbitrary curve intersecting the cocore of the handle at one point. 
A right-handed Dehn twist about a simple closed curve $\alpha \subset S$ is the diffeomorphism $D_\alpha$
that acts on the neighborhood  $N=\alpha \times (0,1)\subset S$ of $\alpha$ as 
$(\theta, t) \mapsto (\theta + 2\pi t, t)$, and fixes  $S \setminus N$, see Figure~\ref{rh-twist}. 
(The term ``positive Dehn twist'' is also common in the literature, but we avoid it since positive Dehn twists 
will correspond to $(-1)$ contact surgeries.)
A left-handed Dehn twist about $\alpha$ is the inverse of $D_\alpha$. 

We recall that the monodromy of an arbitrary open book can be written as a product of left-handed and right-handed 
Dehn twists, and that a contact structure is Stein fillable if and only if it admits an open book with the monodromy 
given by a product of right-handed Dehn twists \cite{Gi}.

\subsection{Contact surgery} \label{cont-surgery}
Let $K$ be a null-homologous Legendrian knot in a contact manifold $(Y, \xi)$. Performing a Dehn surgery on $K$, we cut out a tubular neighborhood of the knot $K$ (i.e. a solid torus) and glue it back in. When the surgery coefficient is $(\pm 1)$ (with respect to the contact framing on $K$), this procedure is compatible with contact structures: the gluing can be done so that the contact structure on the solid torus matches the contact structure on its complement. Moreover, the resulting contact manifold is independent of choices, 
so the $(\pm 1)$ contact surgery is well-defined. Contact surgery is a very useful 
tool, as any contact manifold can be obtained from $(S^3, \xi_{std})$ by a contact surgery on some Legendrian link. We also recall that $(-1)$ contact surgery is in fact the same as Legendrian surgery, while $(+1)$ contact surgery is the operation inverse to it. (Unlike Legendrian surgery, $(+1)$ surgery does not preserve Stein fillability or other similar properties of contact structures.)

Homotopy invariants of a contact structure $\xi$ on $Y$ encode information
on the corresponding plane field. First, we can consider the $\Spinc$ structure $\s$
on $Y$ induced by $\xi$. Secondly, when $c_1(\s)$ is torsion,  
the three-dimensional invariant $d_3(\xi)$ can be defined \cite{Go}. 
 If
$(Y,\xi)$ is the boundary of an almost-complex 4-manifold $(X,J)$,
this invariant is given by
$$
d_3(\xi)=\f{1}{4}(c_1^2(J)-2\chi(X)-3\sign(X)).
$$ 

The homotopy invariants of a contact structure can be read off its contact 
surgery presentation as follows \cite{DGS}.
 
Let $X$ be the four-manifold obtained from $D^4$ by attaching the
2-handles as dictated by the $(\pm 1)$-surgery diagram. 
 Consider an almost-complex structure $J$ defined on
$X$ in the complement of $m$ balls lying in the interior of the
(+1)-handles of $X$. As shown in \cite{DGS}, $J$ induces a
$\Spinc$ structure $s_J$ which extends to all of $X$, and  the $d_3$
invariant of $\xi_L$ can be computed as
\begin{equation}\label{d3}
d_3(\xi_K)=\f{1}{4}(c_1^2(\s_J)-2\chi(X)-3\sign(X))+m.
\end{equation}
This formula is very similar to the case where $(X, J)$ is
almost-complex, except that there is a correction term of $+1$ for
each (+1)-surgery.

Now, suppose that a $2$-handle is attached to the four-manifold
$X$  in the process of Legendrian surgery on a knot $K$, and
denote by $[S]$ the homology class that arises from the Seifert
surface of $K$ capped off inside the handle. It is well-known
\cite{Go} that $c_1(\s_J)$ evaluates on  $[S]$ as the rotation
number of the Legendrian knot $K$. Furthermore, it is shown in
\cite{DGS} that the same result is true for $(+1)$-contact
surgeries (for the $\Spinc$ structure $s_J$ on $X$ described
above).

\subsection{Surgery descriptions from open books} 
There are two ways to describe 
a given contact structure, via an open book decomposition or a contact surgery diagram. 
We will need to switch between the two descriptions.  
A contact surgery diagram consists of 
a Legendrian link in $S^3$ with surgery coefficients. We can find an open book decomposition 
of $S^3$ whose page contains this link. Thus components of the surgery link correspond to 
curves on the page; we perform right-handed Dehn twists on curves 
corresponding to Legendrian surgeries, and left-handed Dehn twists on those corresponding to $(+1)$ 
contact surgeries.  The resulting open book is compatible with the contact structure given by 
the surgery diagram \cite{AO, Pl2, Et3}. Conversely, given an open book decomposition of a given contact manifold, 
we will need to obtain a contact surgery diagram. To this end, we assume that the monodromy 
of the open book contains a sequence of Dehn twists producing $(S^3, \xi_{std})$ (this can always be achieved 
by composing the given monodromy with a few Dehn twists and their inverses). We can then embed the page 
of the open book into $S^3$; if the curves on which the remaining Dehn twists are to be performed become 
Legendrian knots in $S^3$, we can replace the Dehn twists with  $(\pm 1)$ contact surgeries to obtain 
the required surgery diagram.  (Note that a ``compatible'' embedding will imply that the contact framing 
of the Legendrian knots is the same as the page framing.) We perform this procedure in detail in Section~\ref{surg for cov}.

\subsection{The induced contact structure on $\Sigma_p(K)$}\label{subsec-induced-str} 
Given a branched $p$-fold cyclic branched cover $\Sigma_p(K)$ for a transverse knot $K$, we describe the natural contact structure $\xi_p(K)$ on $\Sigma_p(K)$ as follows. In local 
coordinates $(r, \theta, z)$ near the knot $K= \{r=0\}$, we may write 
the covering map as $(r, \theta, z)\mapsto (r^p, p\,\theta, z)$, and assume that 
$\xi = \ker (dz+r^2\,d\theta)$. We set $\xi_p(K)$ to be the kernel of the 
pull-back form; however, the pull-back form $dz+p r^{2p}\, d\theta$ fails to be contact along 
the knot. To avoid this issue, we can define a new contact form by interpolating between the form $dz+r^2 d\theta$ in a small tubular neighborhood of $K$ and the pull-back form 
on the branched cover away from $K$. Its kernel is a contact structure 
which is independent of choices. (This construction is explained in detail in \cite{Pla} for 
branched double covers and works for links and higher order covers with only notational changes.) 

We can also describe the contact structure on  $\Sigma_p(K)$ via open books, by representing $K$ as a braid. We then
consider a branched cover of the standard open book  for $S^3$ whose binding is the braid axis, and page a disk meeting the $n$-braid $K$ at $n$ points. We adopt this approach in the next section, determining how the half-twist generators of the braid $K$ lift to the branched cover. It is clear that the resulting contact structure is isotopic to 
the one described above.

%%%%%%%%%%%%%%%%%%%%%%%%%%%%%%%

\section{Open Books and Surgeries from Braids} \label{obsurg}

\subsection{Dehn Twists and Crossings} 

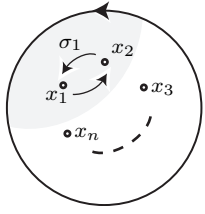
\begin{figure}[htpb!]
\begin{center}
\begin{picture}(75,80)
\put(0,0){\includegraphics{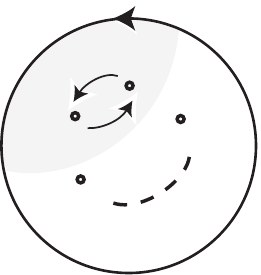}}
\put(15, 40){$\se{x_1}$}  \put(40, 58){$\se{x_2}$}  
\put(26, 26){$\se{x_n}$}  \put(55, 43){$\se{x_3}$}  
\put(20, 60){$\se{\sigma_1}$}
\end{picture}
\caption{A page $D$ and the map $\sigma_1$ acting on the shaded region $U_1$.}\label{sigma}
\end{center}
\end{figure}

Let $K \subset (S^3, \xi_{sym})$ be a transverse link.  Identifying $K$ with a closed braid about the $z$-axis, let $\sigma_{i_1} \sigma_{i_2} \cdots \sigma_{i_K} \in B_n$ denote a braid representation of $K$.  Let $D=\{(r, \theta, z) | \theta=0, r>0 \}\subset \mathbb{R}\cup \{\infty\} = S^3$ be a disk. Then $K$ transversely intersects $D$ in $n$ points $x_1, \cdots, x_n$. We may regard $\sigma_j\in B_n$ as a diffeomorphism of $D$ that exchanges $x_j, x_{j+1}$ in the neighborhood $U_j$ of $x_i, x_{j+1}$ and fixes  $D \setminus U_j$ (Figure~\ref{sigma}). 
Let $\phi_K=\sigma_{i_1} \sigma_{i_2} \cdots \sigma_{i_K}$ be a monodromy map of $D$. The symmetric contact structure $(S^3,\xi_{sym})$ is supported by the open book decomposition $(D, \phi_K)$ of $S^3$, whose binding is the $z$-axis (braid axis) and pages are disks $D$. 

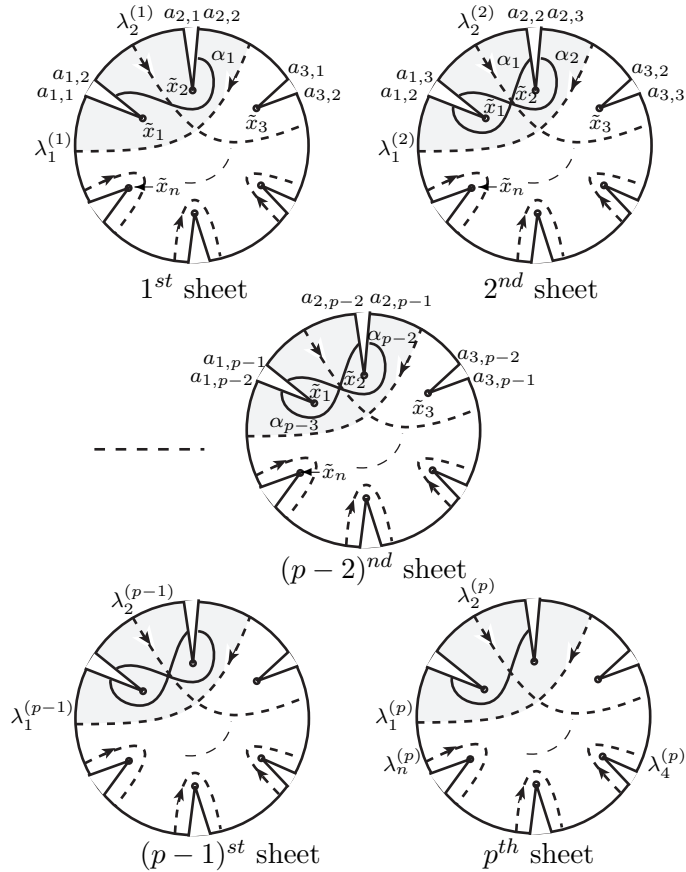
\begin{figure}[htpb!]
\begin{center}
\begin{picture}(218, 325)
\put(0,10){\includegraphics{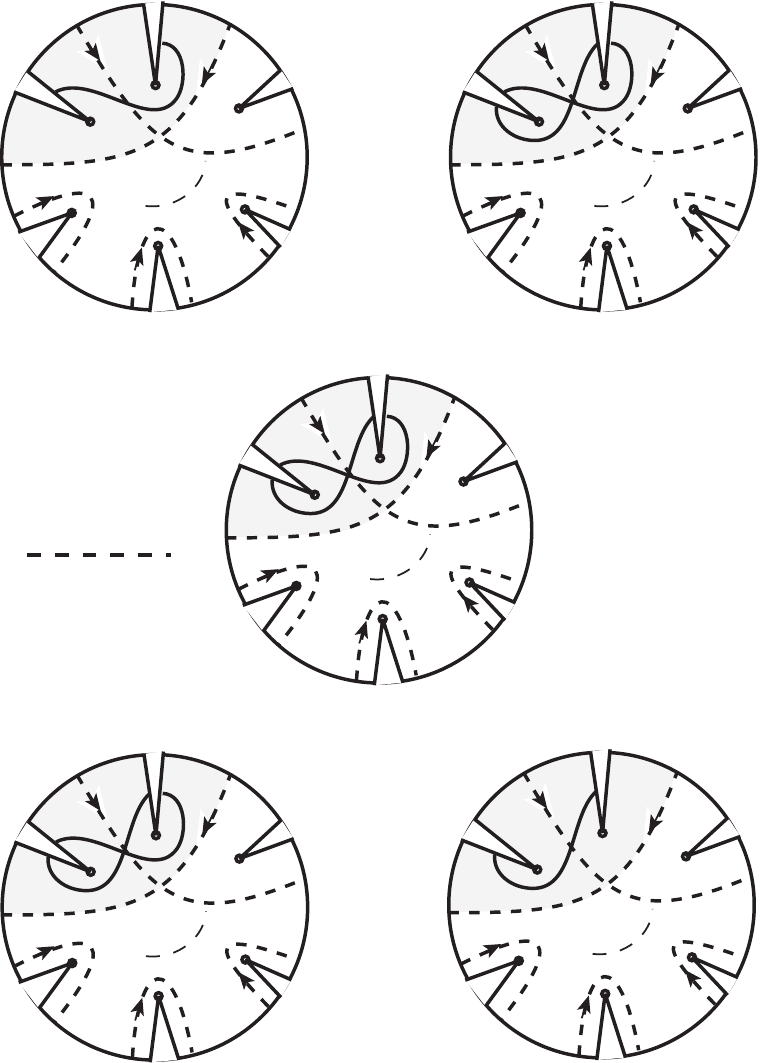}}
% on the way!****
%\psfrag{1}{$\se{\tilde{x}_2}$} \psfrag{2}{$\se{\tilde{x}_1}$} 
%\psfrag{3}{$\se{\tilde{x}_n}$} \psfrag{4}{$\se{\tilde{x}_3}$}  
%\psfrag{5}{$p$-th sheet}        \psfrag{6}{$1$-st sheet} 
%\psfrag{7}{$2$-nd sheet}       \psfrag{8}{$(p-2)$-th sheet} \psfrag{9}{$(p-1)$-th sheet} 
%\psfrag{a}{$\se{\alpha_{p-1}}$}  \psfrag{b}{$\se{\alpha_1}$} 
%\psfrag{c}{$\se{\alpha_2}$}        \psfrag{d}{$\se{\alpha_{p-3}}$} 
%\psfrag{e}{$\se{\alpha_{p-2}}$}  \psfrag{f}{$\se{a_{2,p}}$} 
%\psfrag{g}{$\se{a_{2,1}}$}          \psfrag{h}{$\se{a_{2,2}}$}  
%\psfrag{i}{$\se{a_{2,3}}$}           \psfrag{j}{$\se{a_{2,p-2}}$} 
%\psfrag{k}{$\se{a_{2,p-1}}$}       \psfrag{l}{$\se{a_{1,p}}$} 
%\psfrag{m}{$\se{a_{1,1}}$}         \psfrag{n}{$\se{a_{1,2}}$} 
%\psfrag{o}{$\se{a_{1,3}}$}          \psfrag{p}{$\se{a_{1,p-2}}$} 
%\psfrag{q}{$\se{a_{1,p-1}}$}       \psfrag{r}{$\se{a_{3,p}}$} 
%\psfrag{s}{$\se{a_{3,1}}$}          \psfrag{t}{$\se{a_{3,2}}$} 
%\psfrag{u}{$\se{a_{3,3}}$}          \psfrag{v}{$\se{a_{3,p-2}}$} 
%\psfrag{w}{$\se{a_{3,p-1}}$}
\put(145,99){$\se\lambda^{(p)}_2$}
 \put(114,53){$\se\lambda^{(p)}_1$} 
\put(117, 35){$\se\lambda^{(p)}_n$} 
\put(217, 35){$\se\lambda^{(p)}_{4}$} 
\put(52,304){${\se\alpha_1}$}
\put(34,320){$\se{a_{2,1}}$} 
\put(49,320){$\se{a_{2,2}}$} 
\put(81,299){$\se{a_{3,1}}$} 
\put(87,290){$\se{a_{3,2}}$} 
\put(-14,290){$\se{a_{1,1}}$} 
\put(-8,297){$\se{a_{1,2}}$} 
\put(35,292){$\se{\tilde{x}_2}$} 
\put(26,275){$\se{\tilde{x}_1}$}
\put(65,278){$\se{\tilde{x}_3}$}
\put(31,254){$\se{\tilde{x}_n}$}
\put(30,256){\vector(-1,0){7}}   
\put(25,214){$1^{st}$ sheet}
\put(182,304){${\se\alpha_2}$}
\put(160,302){${\se\alpha_1}$}
\put(164,320){$\se{a_{2,2}}$} 
\put(179,320){$\se{a_{2,3}}$} 
\put(211,299){$\se{a_{3,2}}$} 
\put(217,290){$\se{a_{3,3}}$} 
\put(116,290){$\se{a_{1,2}}$} 
\put(122,297){$\se{a_{1,3}}$} 
\put(167,289){$\se{\tilde{x}_2}$} 
\put(155,284){$\se{\tilde{x}_1}$}
\put(195,278){$\se{\tilde{x}_3}$}
\put(161,254){$\se{\tilde{x}_n}$}
\put(160,256){\vector(-1,0){7}}   
\put(155,214){$2^{nd}$ sheet}
\put(111,198){${\se\alpha_{p-2}}$}
\put(74,165){${\se\alpha_{p-3}}$}
\put(86,212){$\se{a_{2,p-2}}$} 
\put(112,212){$\se{a_{2,p-1}}$} 
\put(145,191){$\se{a_{3,p-2}}$} 
\put(151,182){$\se{a_{3,p-1}}$} 
\put(44,182){$\se{a_{1,p-2}}$} 
\put(49,189){$\se{a_{1,p-1}}$} 
\put(102,181){$\se{\tilde{x}_2}$} 
\put(89,176){$\se{\tilde{x}_1}$}
\put(128,170){$\se{\tilde{x}_3}$}
\put(94,146){$\se{\tilde{x}_n}$}
\put(93,148){\vector(-1,0){7}}   
\put(78,108){$(p-2)^{nd}$ sheet}
\put(16, 316){$\se\lambda^{(1)}_2$} 
\put(-15, 268){$\se\lambda^{(1)}_1$} 
\put(145, 316){$\se\lambda^{(2)}_2$} 
\put(115,268){$\se\lambda^{(2)}_1$} 
\put(14, 97){$\se\lambda^{(p-1)}_2$} 
\put(-24, 53){$\se\lambda^{(p-1)}_1$} 
\put(25,0){$(p-1)^{st}$ sheet}
\put(155,0){$p^{th}$ sheet}
\end{picture}
\caption{A page $\tilde{D}$ and simple closed curves $\alpha_k$'s. The region $\pi^{-1}(U_1)$ is shaded.}\label{Dehn-twist}
\end{center}
\end{figure}

Fix $p\geq 2$, and let $\pi : \Sigma_p(K) \rightarrow S^3$  be the $p$-fold cyclic covering branched along $K$.
The covering $\pi$  induces 
the open book decomposition $(\tilde{D} , \tilde{\phi}_K) = (\pi^{-1}(D), \pi^{-1}(\phi_K))$ of 
$\Sigma_p(K)$ given by the lift of the open book $(D, \phi_K)$. 
The surface $\tilde{D}$ can be obtained by gluing  $p$ copies of $D$ along slits as in Figure~\ref{Dehn-twist}. 
For example, denoting $\tilde{x}_j = \pi^{-1}(x_j) \in \tilde{D}$, we identify the edge $a_{j,k} \tilde{x}_j$ of the $k$-th sheet with the edge $a_{j,k} \tilde{x}_j$ of the $(k+1)$-th sheet. To compute the monodromy map $\tilde{\phi}_K$ we need the following lemma.

%%%%% Dehn twist lemma  %%%%%%%%%%%%%%
\begin{lemma}\label{Dehn twist lemma}

Let $\alpha_k \subset \tilde{D}$ for $k=1,\cdots, p-1$ be a simple closed curve as in Figure \ref{Dehn-twist}. Let 
$D_k= D_{\alpha_k}$ be the right-handed Dehn twist along $\alpha_k$. 
Then the lift $\tilde{\sigma_1}$ of $\sigma_1$ is $D_1 \circ D_2 \circ \cdots \circ D_{p-1}$ 
(where  $D_{p-1}$ comes first and  $D_1$ last). 
\end{lemma}

\begin{proof}
For simplicity, denote $\sigma := \sigma_1$ and $U:= U_1$.
We need to show  that up to isotopy,   
\begin{equation}\label{eq}
\pi \circ \sigma^{-1} \circ \pi \circ D_1 \circ D_2 \circ \cdots \circ D_{p-1} = id_{\tilde{D}}.
\end{equation}

\begin{figure}[htpb!]
\begin{center}
\begin{picture}(334,430)
\put(0,0){\includegraphics{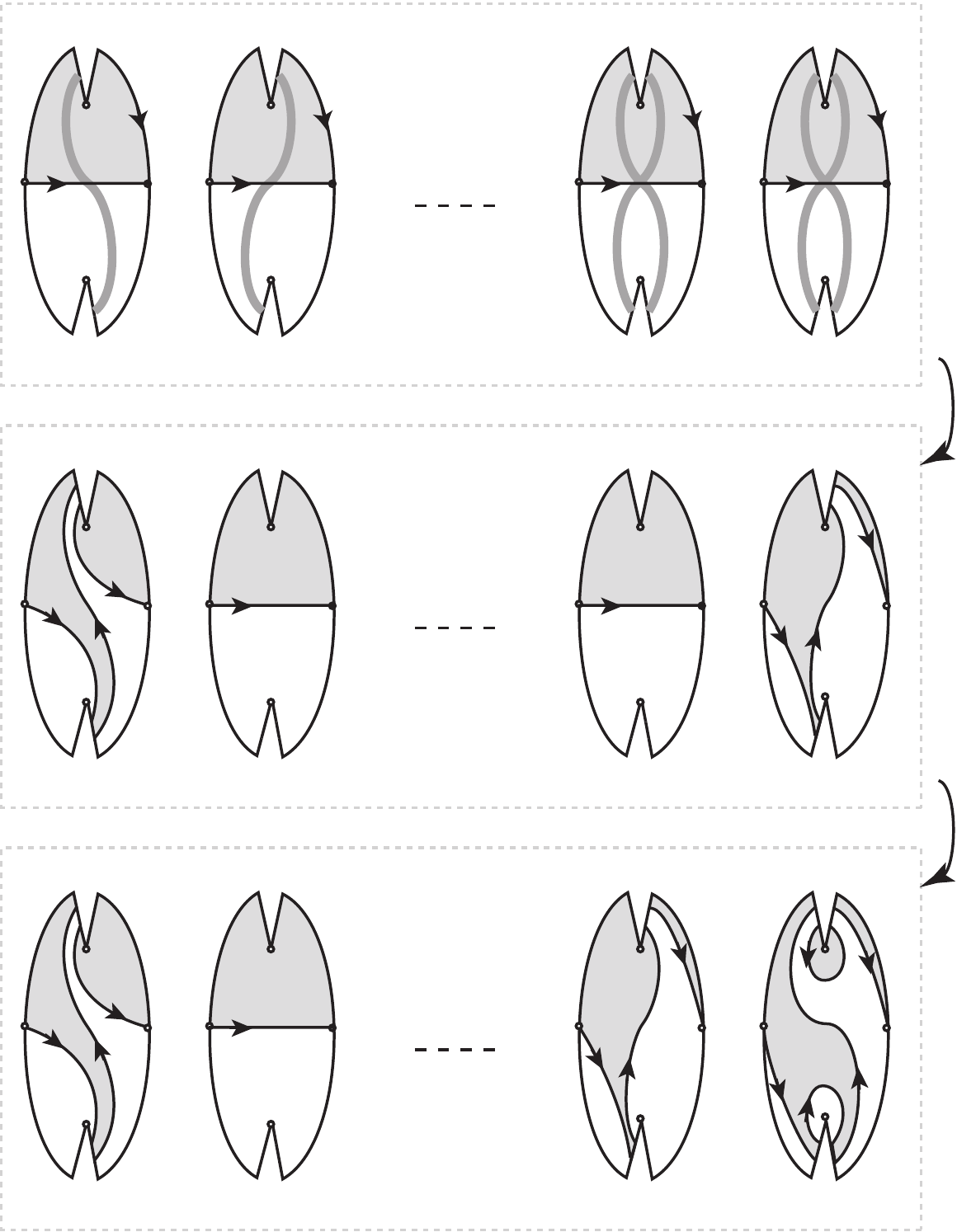}}
%top line
\put(14,418){$\se{\alpha_{2,p}}$}
\put(32,418){$\se{\alpha_{2,1}}$}
\put(78,418){$\se{\alpha_{2,1}}$}
\put(96,418){$\se{\alpha_{2,2}}$}
\put(201,418){$\se{\alpha_{2,p-2}}$}
\put(228,418){$\se{\alpha_{2,p-1}}$}
\put(266,418){$\se{\alpha_{2,p-1}}$}
\put(293,418){$\se{\alpha_{2,p-2}}$}
%next line
\put(14,306){$\se{\alpha_{1,1}}$}
\put(32,306){$\se{\alpha_{2,1}}$}
\put(78,306){$\se{\alpha_{2,1}}$}
\put(96,306){$\se{\alpha_{1,1}}$}
\put(201,306){$\se{\alpha_{1,p-1}}$}
\put(228,306){$\se{\alpha_{1,p-2}}$}
\put(266,306){$\se{\alpha_{1,p-2}}$}
\put(293,306){$\se{\alpha_{1,1}}$}
\put(337,288){$D_{p-1}$}
\put(14,272){$\se{\alpha_{2,p}}$}
\put(32,272){$\se{\alpha_{2,1}}$}
\put(78,272){$\se{\alpha_{2,1}}$}
\put(96,272){$\se{\alpha_{2,2}}$}
\put(201,272){$\se{\alpha_{2,p-2}}$}
\put(228,272){$\se{\alpha_{2,p-1}}$}
\put(266,272){$\se{\alpha_{2,p-1}}$}
\put(293,272){$\se{\alpha_{2,p-2}}$}
\put(14,160){$\se{\alpha_{1,1}}$}
\put(32,160){$\se{\alpha_{2,1}}$}
\put(78,160){$\se{\alpha_{2,1}}$}
\put(96,160){$\se{\alpha_{1,1}}$}
\put(201,160){$\se{\alpha_{1,p-1}}$}
\put(228,160){$\se{\alpha_{1,p-2}}$}
\put(266,160){$\se{\alpha_{1,p-2}}$}
\put(293,160){$\se{\alpha_{1,1}}$}
\put(337,141){$D_{p-2}$}
\put(14,124){$\se{\alpha_{2,p}}$}
\put(32,124){$\se{\alpha_{2,1}}$}
\put(78,124){$\se{\alpha_{2,1}}$}
\put(96,124){$\se{\alpha_{2,2}}$}
\put(201,124){$\se{\alpha_{2,p-2}}$}
\put(228,124){$\se{\alpha_{2,p-1}}$}
\put(266,124){$\se{\alpha_{2,p-1}}$}
\put(293,124){$\se{\alpha_{2,p-2}}$}
\put(14,12){$\se{\alpha_{1,1}}$}
\put(32,12){$\se{\alpha_{2,1}}$}
\put(78,12){$\se{\alpha_{2,1}}$}
\put(96,12){$\se{\alpha_{1,1}}$}
\put(201,12){$\se{\alpha_{1,p-1}}$}
\put(228,12){$\se{\alpha_{1,p-2}}$}
\put(266,12){$\se{\alpha_{1,p-2}}$}
\put(293,12){$\se{\alpha_{1,1}}$}
%x_2
\put(28,386){$\se{\tilde{x}_2}$}
\put(93,386){$\se{\tilde{x}_2}$}
\put(221,386){$\se{\tilde{x}_2}$}
\put(286,386){$\se{\tilde{x}_2}$}
\put(28,239){$\se{\tilde{x}_2}$}
\put(93,239){$\se{\tilde{x}_2}$}
\put(221,239){$\se{\tilde{x}_2}$}
\put(286,239){$\se{\tilde{x}_2}$}
\put(28,91){$\se{\tilde{x}_2}$}
\put(93,91){$\se{\tilde{x}_2}$}
\put(221,91){$\se{\tilde{x}_2}$}
\put(286,91){$\se{\tilde{x}_2}$}
% X_1
\put(28,336){$\se{\tilde{x}_1}$}
\put(93,336){$\se{\tilde{x}_1}$}
\put(221,336){$\se{\tilde{x}_1}$}
\put(286,336){$\se{\tilde{x}_1}$}
\put(25,189){$\se{\tilde{x}_1}$}
\put(93,189){$\se{\tilde{x}_1}$}
\put(221,189){$\se{\tilde{x}_1}$}
\put(288,190){$\se{\tilde{x}_1}$}
\put(25,42){$\se{\tilde{x}_1}$}
\put(93,42){$\se{\tilde{x}_1}$}
\put(223,42){$\se{\tilde{x}_1}$}
\put(286,43){$\se{\tilde{x}_1}$}
%lambdas
\put(11,373){$\se{\lambda_2^{(p)}}$}
\put(76,373){$\se{\lambda_2^{(1)}}$}
\put(270,373){$\se{\lambda_2^{(p-1)}}$}
\end{picture}
%\psfrag{1}{$\se{\tilde{x}_2}$}      \psfrag{2}{$\se{\tilde{x}_1}$} 
%\psfrag{a}{$\se{\alpha_{p-1}}$}   \psfrag{b}{$\se{\alpha_1}$} 
%\psfrag{d}{$\se{\alpha_{p-3}}$}   \psfrag{e}{$\se{\alpha_{p-2}}$} 
%\psfrag{D}{$\se{D_{p-1}}$}          \psfrag{E}{$\se{D_{p-2}}$} 
%\psfrag{f}{$\se{a_{2,p}}$}            \psfrag{g}{$\se{a_{2,1}}$} 
%\psfrag{h}{$\se{a_{2,2}}$}           \psfrag{i}{$\se{a_{2,3}}$} 
%\psfrag{j}{$\se{a_{2,p-2}}$}         \psfrag{k}{$\se{a_{2,p-1}}$}
%\psfrag{l}{$\se{a_{1,1}}$}            \psfrag{m}{$\se{a_{1,p}}$} 
%\psfrag{n}{$\se{a_{1,1}}$}           \psfrag{p}{$\se{a_{1,p-1}}$} 
%\psfrag{q}{$\se{a_{1,p-2}}$}
%\psfrag{A}{$\se\lambda^{(p)}_2$} \psfrag{B}{$\se\lambda^{(p)}_1$} 
%\psfrag{G}{$\se\lambda^{(1)}_2$} \psfrag{H}{$\se\lambda^{(1)}_1$}
%\psfrag{T}{$\se\lambda^{(p-1)}_2$} \psfrag{Q}{$\se\lambda^{(p-1)}_1$}
%\includegraphics[height=15cm]{alpha-p}
\caption{Actions of $D_{p-1}$ and $D_{p-2}$ on $\pi^{-1}(U)$. }\label{alpha-p}
\end{center}
\end{figure}
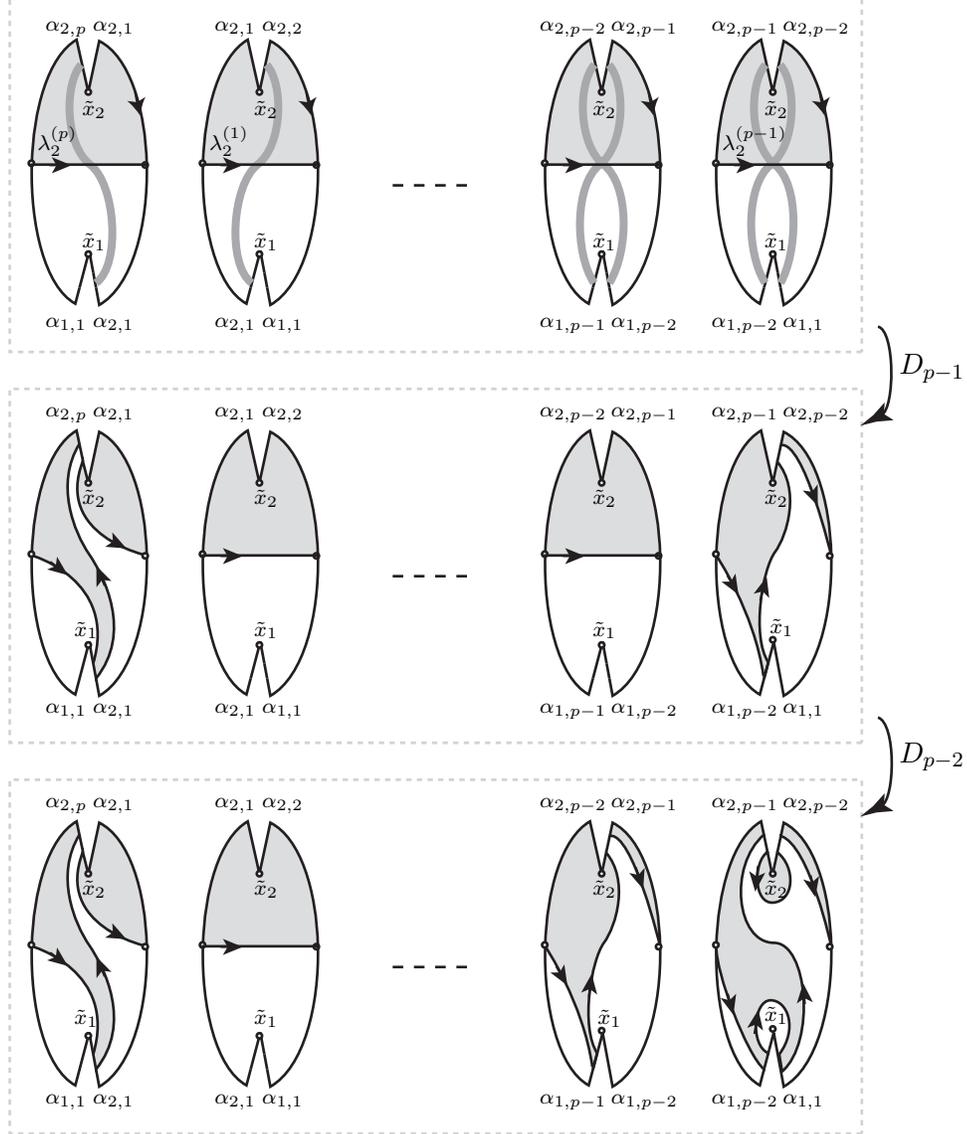
Cut $\tilde{D}$ into $n+p$ disks along oriented properly embedded arcs $\lambda^{(i)}_j$ where $i=1, \cdots, p$ and $j=1,\cdots, n-1$, dashed in Figure~\ref{Dehn-twist}. We will check that after an isotopy, 
the map $\pi \circ \sigma^{-1} \circ \pi \circ D_1 \circ D_2 \circ \cdots \circ D_{p-1}$ fixes each vertex and oriented edge of the graph $\cup_{i,j} \lambda^{(i)}_j$.
Our statement will then follow from the  Alexander method (see e.g. \cite[Proposition 3.4]{FM}). The Alexander method is based on the observation that a diffeomorphism of $D$ fixing $\partial D$ is isotopic to identity; this observation is applied to each of the $(n+p)$ disks.

Since the Dehn twists are performed on curves $\alpha_1, \dots, \alpha_{p-1}$ which all lie 
 in $\pi^{-1}(U)$,  we can assume that all the $\lambda$-arcs except 
 $\lambda^{(i)}_1$'s are fixed by $D_1 \circ D_2 \circ \cdots \circ D_{p-1}$. 
\begin{figure}[htpb!]
\begin{center}
\begin{picture}(366,120)
%\psfrag{1}{$\se{\tilde{x}_2}$} \psfrag{2}{$\se{\tilde{x}_1}$} 
%\psfrag{f}{$\se{a_{2,p}}$} \psfrag{g}{$\se{a_{2,1}}$} 
%\psfrag{h}{$\se{a_{2,2}}$} \psfrag{i}{$\se{a_{2,3}}$} 
%\psfrag{j}{$\se{a_{2,p-2}}$} \psfrag{k}{$\se{a_{2,p-1}}$}
%\psfrag{l}{$\se{a_{1,1}}$}            \psfrag{m}{$\se{a_{1,p}}$} 
%\psfrag{n}{$\se{a_{1,2}}$}           \psfrag{o}{$\se{a_{1,3}}$} 
%\psfrag{p}{$\se{a_{1,p-1}}$}         \psfrag{q}{$\se{a_{1,p-2}}$}
%\put(0,0){\includegraphics[height=4cm]{pre-image2}}
\put(0,10){\includegraphics{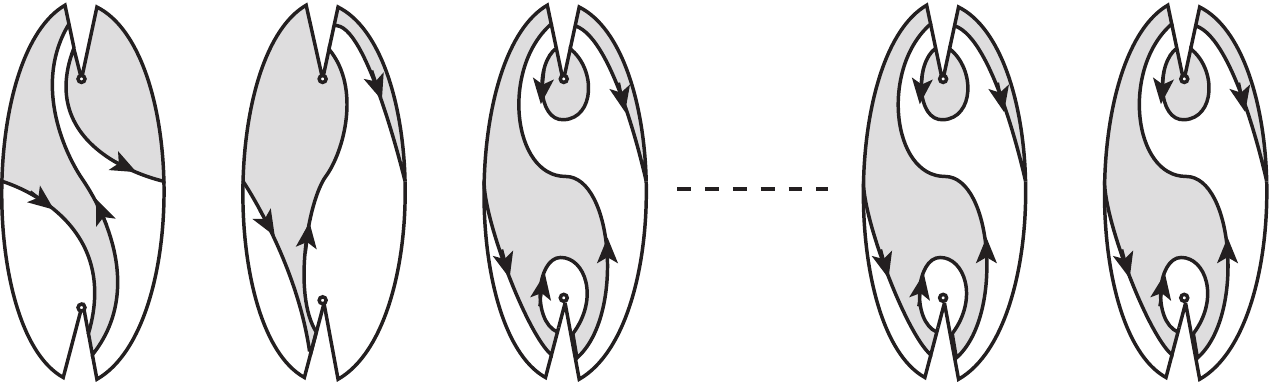}}
\put(8,122){$\se{\alpha_{2,p}}$}
\put(26,122){$\se{\alpha_{2,1}}$}
\put(78,122){$\se{\alpha_{2,1}}$}
\put(96,122){$\se{\alpha_{2,2}}$}
\put(148,122){$\se{\alpha_{2,2}}$}
\put(166,122){$\se{\alpha_{2,3}}$}
\put(248,122){$\se{\alpha_{2,p-2}}$}
\put(275,122){$\se{\alpha_{2,p-1}}$}
\put(320,122){$\se{\alpha_{2,p-1}}$}
\put(346,122){$\se{\alpha_{2,p}}$}
\put(8,4){$\se{\alpha_{1,1}}$}
\put(26,4){$\se{\alpha_{1,p}}$}
\put(78,4){$\se{\alpha_{1,2}}$}
\put(96,4){$\se{\alpha_{1,1}}$}
\put(148,4){$\se{\alpha_{1,3}}$}
\put(166,4){$\se{\alpha_{1,2}}$}
\put(248,4){$\se{\alpha_{1,p-1}}$}
\put(275,4){$\se{\alpha_{1,p-2}}$}
\put(324,4){$\se{\alpha_{1,p}}$}
\put(346,4){$\se{\alpha_{1,p-1}}$}
\put(21,90){$\se{\tilde{x}_2}$}
\put(90,90){$\se{\tilde{x}_2}$}
\put(159,90){$\se{\tilde{x}_2}$}
\put(269,90){$\se{\tilde{x}_2}$}
\put(338,90){$\se{\tilde{x}_2}$}
\put(19,36){$\se{\tilde{x}_1}$}
\put(90,36){$\se{\tilde{x}_1}$}
\put(159,37){$\se{\tilde{x}_1}$}
\put(269,37){$\se{\tilde{x}_1}$}
\put(338,37){$\se{\tilde{x}_1}$}
\end{picture}
\caption{The region $W =D_1 \circ D_2 \circ \cdots \circ D_{p-1}(\pi^{-1}(U))$. }\label{pre-image2}
\end{center}
\end{figure}
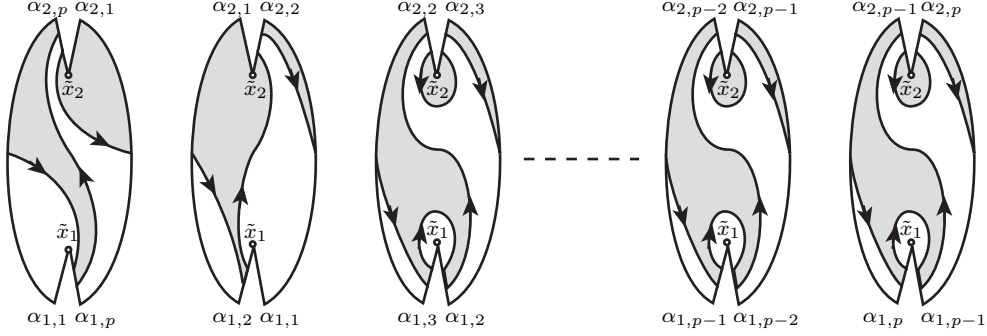
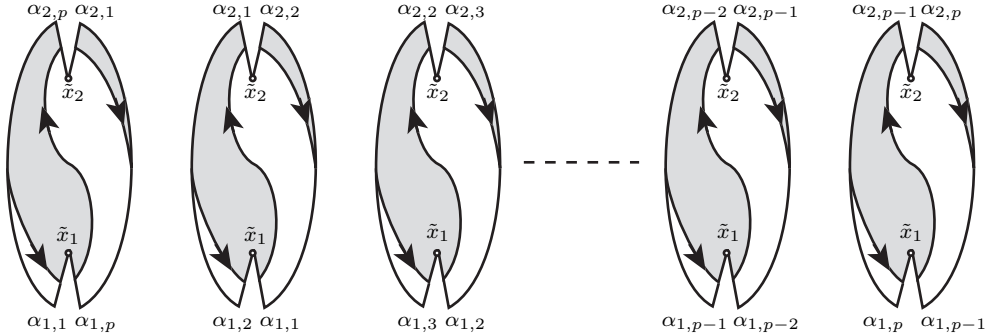
\begin{figure}[htpb!]
\begin{center}
\begin{picture}(366,120)
\put(0,10){\includegraphics{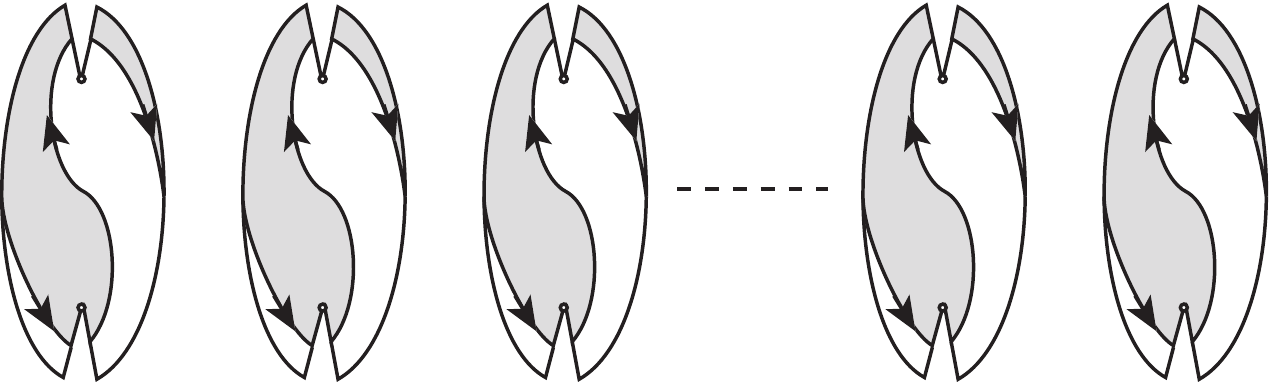}}
\put(8,122){$\se{\alpha_{2,p}}$}
\put(26,122){$\se{\alpha_{2,1}}$}
\put(78,122){$\se{\alpha_{2,1}}$}
\put(96,122){$\se{\alpha_{2,2}}$}
\put(148,122){$\se{\alpha_{2,2}}$}
\put(166,122){$\se{\alpha_{2,3}}$}
\put(248,122){$\se{\alpha_{2,p-2}}$}
\put(275,122){$\se{\alpha_{2,p-1}}$}
\put(320,122){$\se{\alpha_{2,p-1}}$}
\put(346,122){$\se{\alpha_{2,p}}$}
\put(8,4){$\se{\alpha_{1,1}}$}
\put(26,4){$\se{\alpha_{1,p}}$}
\put(78,4){$\se{\alpha_{1,2}}$}
\put(96,4){$\se{\alpha_{1,1}}$}
\put(148,4){$\se{\alpha_{1,3}}$}
\put(166,4){$\se{\alpha_{1,2}}$}
\put(248,4){$\se{\alpha_{1,p-1}}$}
\put(275,4){$\se{\alpha_{1,p-2}}$}
\put(324,4){$\se{\alpha_{1,p}}$}
\put(346,4){$\se{\alpha_{1,p-1}}$}
\put(21,90){$\se{\tilde{x}_2}$}
\put(90,90){$\se{\tilde{x}_2}$}
\put(159,90){$\se{\tilde{x}_2}$}
\put(269,90){$\se{\tilde{x}_2}$}
\put(338,90){$\se{\tilde{x}_2}$}
\put(19,36){$\se{\tilde{x}_1}$}
\put(90,36){$\se{\tilde{x}_1}$}
\put(159,37){$\se{\tilde{x}_1}$}
\put(269,37){$\se{\tilde{x}_1}$}
\put(338,37){$\se{\tilde{x}_1}$}
\end{picture}
%
%\psfrag{1}{$\se{\tilde{x}_2}$}   \psfrag{2}{$\se{\tilde{x}_1}$} 
%\psfrag{f}{$\se{a_{2,p}}$}         \psfrag{g}{$\se{a_{2,1}}$} 
%\psfrag{h}{$\se{a_{2,2}}$}         \psfrag{i}{$\se{a_{2,3}}$} 
%\psfrag{j}{$\se{a_{2,p-2}}$}        \psfrag{k}{$\se{a_{2,p-1}}$}
%\psfrag{l}{$\se{a_{1,1}}$}            \psfrag{m}{$\se{a_{1,p}}$} 
%\psfrag{n}{$\se{a_{1,2}}$}           \psfrag{o}{$\se{a_{1,3}}$} 
%\psfrag{p}{$\se{a_{1,p-1}}$}         \psfrag{q}{$\se{a_{1,p-2}}$}
%\includegraphics[height=4cm]{pre-image}
\caption{The region $W'$ obtained by finger moves applied to $W$.}\label{pre-image}
\end{center}
\end{figure}

Therefore, we focus on $\pi^{-1}(U)$ shown  in the top box of Figure~\ref{alpha-p} to understand how the arcs $\lambda^{(i)}_1$ change under the map $\pi \circ \sigma^{-1} \circ \pi \circ D_1 \circ D_2 \circ \cdots \circ D_{p-1}$. (Note that $\pi^{-1}(U)$ is the union of two shaded regions in Figure~\ref{Dehn-twist}, while in Figure~\ref{alpha-p} one of these regions is shaded and the other is white; the two regions are separated by the arcs  $\lambda^{(i)}_1$). The Dehn twist $D_{p-1}$ changes the region $\pi^{-1}(U)$, shown in the top box of Figure~\ref{alpha-p}, to $D_{p-1}( \pi^{-1}(U))$  as in the second box. Then $D_{p-2}$ changes it to $D_{p-2} (D_{p-1}( \pi^{-1}(U)))$ as in the bottom box. Applying all the Dehn twists $D_1 \circ D_2 \circ \cdots \circ D_{p-1}$ to $\pi^{-1}(U)$, we obtain the region $W$ shown in Figure~\ref{pre-image2}. Next we isotope $W$ fixing the boundary of $W$ by a combination of two local finger moves near $\tilde{x}_1$ and $\tilde{x}_2$, and obtain a region $W'$ as in Figure~\ref{pre-image}.

To complete the proof, we observe that the region $W'$ is precisely $\pi^{-1}( \sigma(U))$.
\end{proof}

Applying the above lemma repeatedly for different pairs of points $x_j, x_{j+1}$, we can write down the monodromy of an arbitrary braid.  We denote the curve $\alpha_k$  introduced in Lemma~\ref{Dehn twist lemma} by $\alpha_k^{j}$ ($k=1, \cdots, p-1$, $j=1, \cdots, n-1)$ when it is related to the twist of branch points $x_j, x_{j+1}$ and lies on the $k$-th and $(k+1)$-th sheets, and write $D_k^j$ for the right-handed Dehn twist around $\alpha_k^{j}$.   
In particular, the $\alpha_k$ curve in Figure~\ref{Dehn-twist} is renamed as $\alpha_k^{1}$, and the corresponding 
Dehn twist is  $D_k^1$.

\begin{prop} \label{sigma1} 
Let $K$ be the braid $\sigma_1 \sigma_2 \dots \sigma_{n-1}$. 
Then the open book for the $p$-fold cover of $K$ given by  Lemma \ref{Dehn twist lemma} is the same as the open book 
for $S^3$ induced by the $(n, p)$-torus link fibration; moreover, the images of the curves $\alpha^{j}_k$
on the Seifert surface of this torus link are as shown on Figure \ref{torus-link}. Each $\alpha^{j}_k$
is an unknot in $S^3$, with page framing $=-1$.  
\end{prop}

\begin{proof}
We first observe that the branched $p$-fold cover of $K$ is $(S^3, \xi_{std})$. This is easy to see: since $K$ is 
the transverse   unknot with $sl=-1$, it can be thought of as the binding of an open book
decomposition of $S^3$ whose page is a disk. The branched $p$-fold cover, then, is given by 
the same open book for any $p$, yielding the standard contact structure on $S^3$.

Lemma \ref{Dehn twist lemma} produces a different open book for the $p$-fold cover of $K$.
The page of this open book, together with the curves $\alpha^{j}_k$,  
can be embedded into $S^3$ as a Seifert surface of the $(p, n)$-torus link shown on Figure \ref{torus-link}.
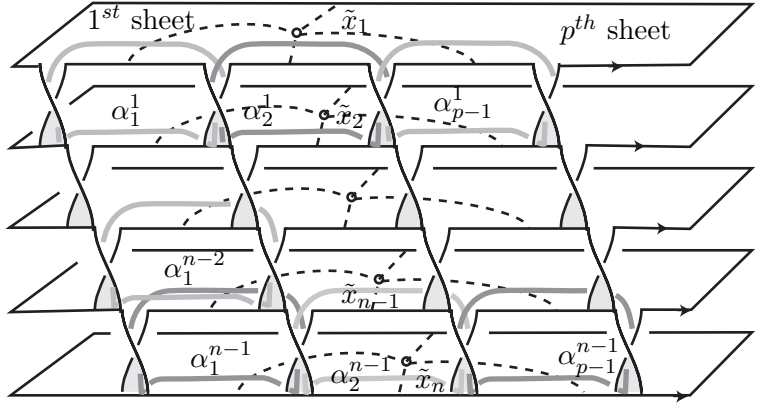
\begin{figure}[htpb!]
\begin{center}
\begin{picture}(288, 160)
\put(0,0){\includegraphics{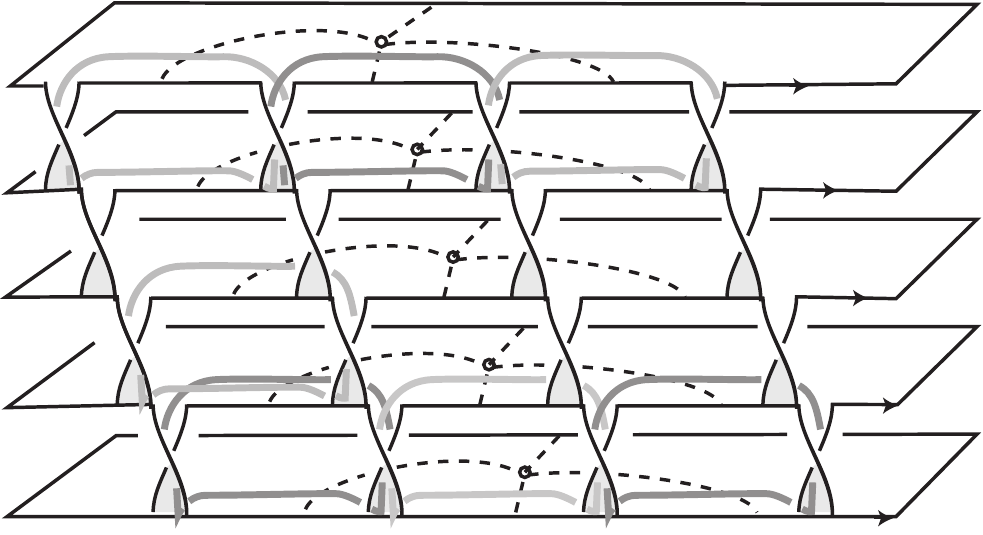}}
\put(155, 8){$\tilde{x}_n$} 
\put(127, 40){$\tilde{x}_{n-1}$} 
\put(125, 108){$\tilde{x}_2$} 
\put(127, 144){$\tilde{x}_1$} 
\put(70, 16){$\alpha_1^{n-1}$} 
\put(123, 10){$\alpha_2^{n-1}$} 
\put(209, 17){$\alpha_{p-1}^{n-1}$} 
\put(60, 50){$\alpha_1^{n-2}$}
\put(40, 110){$\alpha_1^{1}$} 
\put(90, 110){$\alpha_2^{1}$} 
\put(162, 113){$\alpha_{p-1}^{1}$} 
\put(210, 140){$p^{th}$ sheet} 
\put(31, 144){$1^{st}$ sheet} 
\end{picture}
\caption{A page $\tilde{D}$ of the open book with  simple closed curves $\alpha_k^j$ where $p=4$ and $n=5$.} \label{torus-link}
\end{center}
\end{figure}   
It is then clear that each $\alpha^j_k$ is an unknot, with page framing =$-1$. We claim that the torus knot fibration induces the monodromy of the open book given by Lemma \ref{Dehn twist lemma}, i.e. 
 the monodromy of this torus knot is the product of the Dehn twists  
$(D^{n-1}_1 \circ \cdots \circ D^{n-1}_{p-1}) \circ \cdots \circ (D^2_1 \circ \cdots \circ D^2_{p-1})\circ (D^{1}_1 \circ \cdots \circ D^1_{p-1})$.
Since the fiber surface of the torus knot can be obtained by plumbing together a sequence of right-handed Hopf bands whose core curves are  $\alpha_k^{j}$, it is clear that the monodromy of the torus knot is given by a composition 
of the right-handed Dehn twists $D_k^j$. We need to determine the order in which the Dehn twists are performed.

\begin{figure}[htpb!] \begin{center}
\includegraphics[scale=0.70]{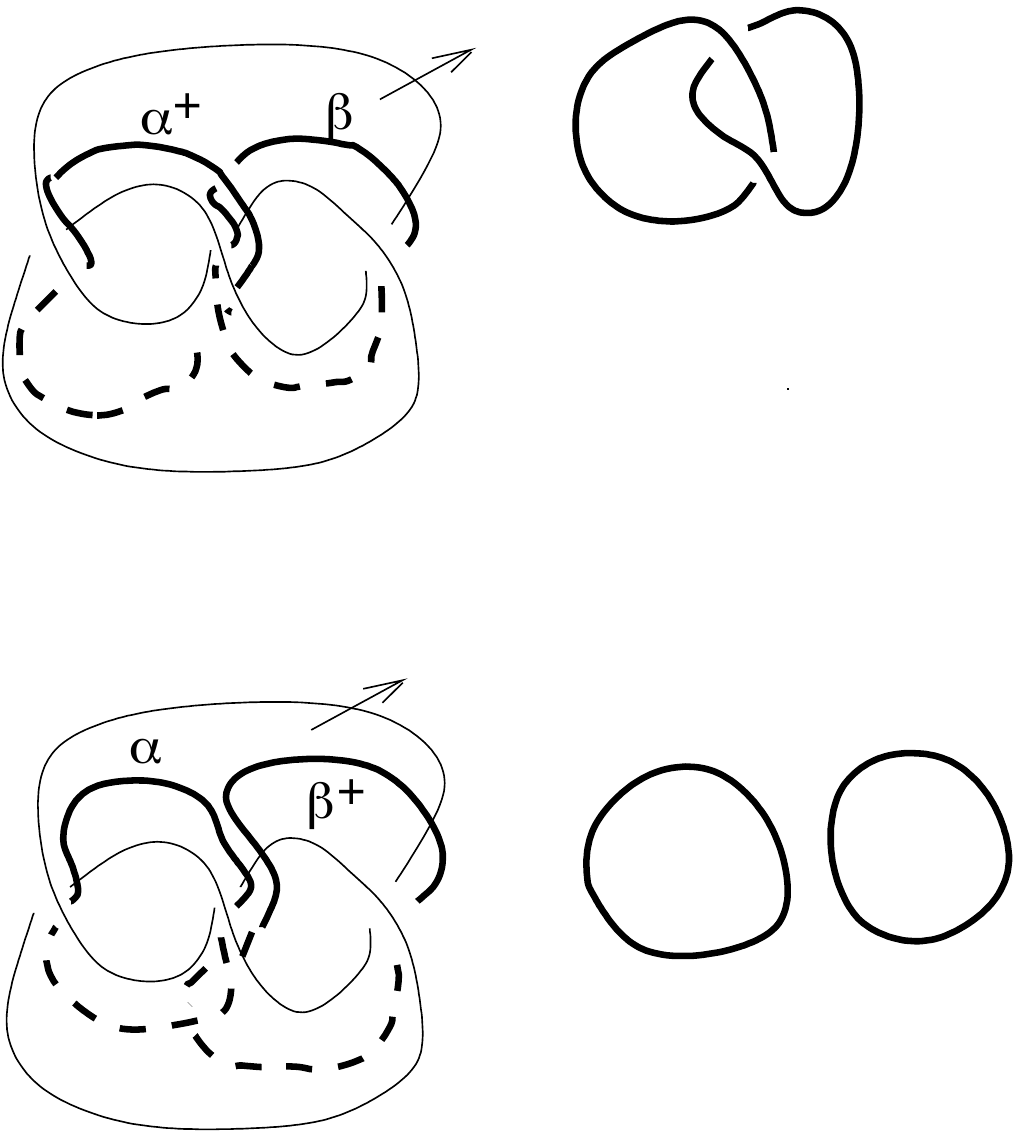}
\caption{The curves $\alpha$ and $\beta$ lie on a fiber of the trefoil knot fibration; $\alpha^+$ and 
$\beta^+$ are their push-offs in the positive normal direction.}\label{linking}
\end{center}\end{figure}  

To simplify the picture, we consider a model example where $n=2$, $p=3$. 
Let $T$ the right-handed trefoil knot and consider the fibration $S^3 \setminus T \to S^1$.
Its monodromy is the product of the Dehn twists around the  curves 
$\alpha = \alpha^1_1$ and $\beta= \alpha^1_2$. Let $P_\theta$, $\theta \in [0, 2\pi)$ be  pages of the corresponding 
open book.  Assume that the curves $\alpha$ and $\beta$ both lie on  $P_{0}$; let $\alpha^+$ and $\beta^+$ 
be their push-offs to the page $P_{\theta^+}$ for some small $\theta^+>0$. Since $S^3 \setminus T$ is oriented 
as a mapping torus, this means that the curves are pushed off in the direction shown by arrow. 
Observe that $\alpha^+$ and $\beta$ form a Hopf link, while $\alpha$ and $\beta^+$ are not linked. 
Suppose that  the monodromy of the pictured trefoil is $D_\beta \circ D_\alpha$, 
and compose it with $D_\alpha^{-1} \circ D_\beta^{-1}$. 
The result is of course the open book with trivial monodromy, which 
gives ${\#_2} S^1 \times S^2$. On the other hand, 
the composition of the two additional Dehn twists corresponds to an integral surgery on $S^3$
performed on the link $\alpha^+ \cup \beta$ (since $D_\alpha^{-1}$ follows $D_\beta^{-1}$,
we need to place a copy of $\alpha$ on the page following the page with $\beta$).
The surgery coefficients are given by (page framing)+1, so we perform 0-surgery on both $\alpha^+$ and $\beta$;
but this surgery on the Hopf link produces $S^3$, not ${\#_2} S^1 \times S^2$. By contrast, 
if we perform 0-surgeries on $\alpha$ and $\beta^+$ which form a trivial link (and correspond to composing 
the trefoil monodromy with $D_\beta^{-1} \circ D_\alpha^{-1}$), we obtain 
${\#_2} S^1 \times S^2$, so we conclude that the trefoil monodromy is $D_\alpha \circ D_\beta$.

Similar argument for various pairs of curves $\alpha^{j}_k$ shows that the monodromy of the torus knot on Figure \ref{torus-link}
is indeed  $(D^{n-1}_1 \circ \cdots \circ D^{n-1}_{p-1}) \circ \cdots \circ (D^2_1 \circ \cdots \circ D^2_{p-1}) \circ (D^{1}_1 \circ \cdots \circ D^1_{p-1})$. The curves  $\alpha^{j}_k$ and ${\alpha^{i}_l}^+$, the push-off 
of  $\alpha^{i}_l$, form a Hopf link whenever 
$(i,l)=(j,k), (j-1,k-1), (j, k-1)$ or $(j-1, k)$, and the trivial unlink otherwise. 
(We return to this in Remark \ref{rmk3.3}, 
see Figure~\ref{Hopf-links} for details.)
\end{proof}

%%%%%%%%%%%%%%%%%%%%%%%%

\subsection{Surgery diagrams for branched covers} \label{surg for cov}

Open books from the previous section will allow us to construct contact surgery diagrams
for the branched covers. In Proposition~\ref{sigma1}, we saw that the branched $p$-fold cover for the transverse  braid $K=\sigma_1 \sigma_2 \dots \sigma_{n-1}$ is $(S^3, \xi_{std})$. Now
consider a transverse $n$-braid $L =\sigma_1 \sigma_2 \dots \sigma_{n-1} b$, where $b$ is an arbitrary braid word. 
The branched $p$-fold cover for $L$ can be obtained from the branched cover for $K$ by performing additional Dehn twists about curves $\alpha_j^i$ in the open book decomposition of $(S^3, \xi_{std})$ considered in Lemma~\ref{Dehn twist lemma}.

The goal of this subsection is to interpret these Dehn twists as contact surgeries. Forgetting the contact structure, we can translate Dehn twists into Dehn surgeries along push-offs of the curves $\alpha_j^i$ to successive pages of our open book. We perform $0$-surgeries for left-handed and $(-2)$-surgeries for right-handed Dehn twists. The order of push-offs is determined by the order of Dehn twists, which in turn is dictated by the braid word $b$ and Lemma~ \ref{Dehn twist lemma}. 

\begin{figure}[ht] 
\begin{center}
\begin{picture}(350, 215)
\put(0,0){\includegraphics{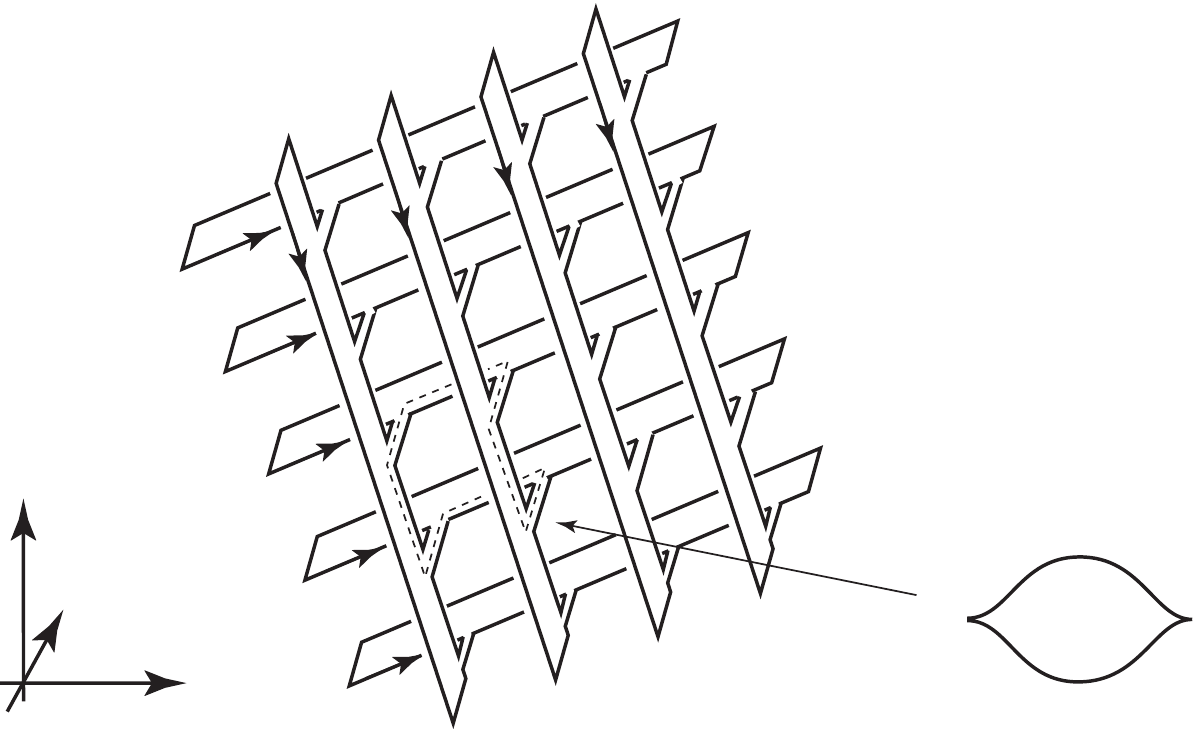}}
\put(55, 15){$x$}
\put(20, 42){$y$}
\put(-8, 62){$z$}
\put(117, 78){$\alpha_1^{n-2}$}
\end{picture}
\end{center}
\caption{A Legendrian realization of Figure~\ref{torus-link}.}\label{LegReal} 
\end{figure}
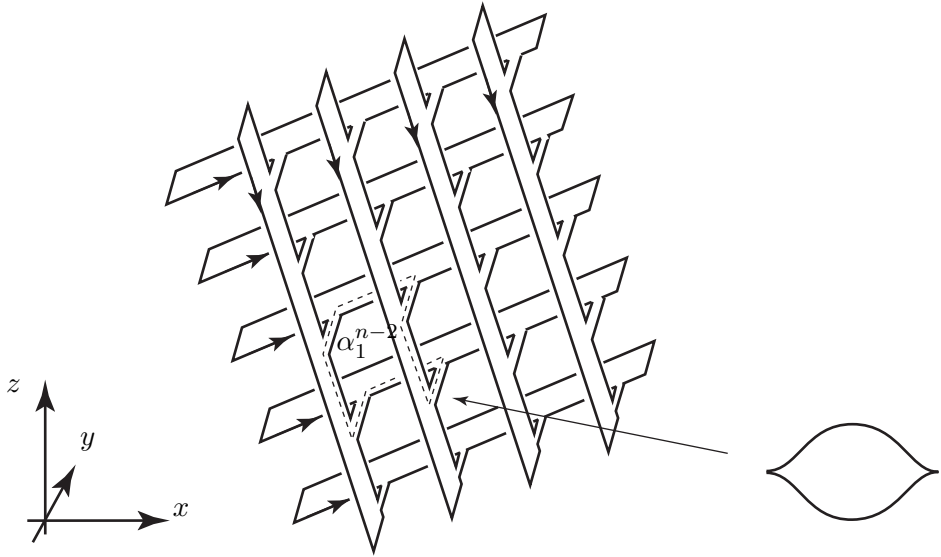

Using Honda's Legendrian Realization \cite{Ho}, we can in principle find an isotopy that  takes all  $\alpha_j^i$ to Legendrian curves whose contact framing matches the page framing, so that $0$- and $(-2)$-surgeries become contact $(\pm 1)$-surgeries. This is almost what we need, but we want an explicit surgery diagram; to this end, we give an explicit Legendrian realization of our curves. Indeed, following \cite{AO} (see \cite[Appendix]{Pl2} for the same construction in the presence of a contact structure), we can embed the fiber surface of a torus link (Figure~\ref{torus-link}) into $S^3$ as the  page $P_0$ of an open book decomposition compatible with $\xi_{std}$, and such that $\alpha_j^i$ are all Legendrian unknots with $tb=-1$. We simply draw this surface as in Figure \ref{LegReal} (assuming as usual that $\xi_{std} = \ker (dz-ydx))$. Various Legendrian push-offs of $\alpha_j^i$ can then be thought of as lying on different pages of the same open book.

To produce a contact surgery diagram of the $p$-fold branched cover for a transverse braid $L=(\sigma_1 \sigma_2 \dots \sigma_{n-1}) b $, we can now start with $S^3$, write down the monodromy of the open book as dictated by the crossings of $b$ and Lemma~\ref{Dehn twist lemma}, and then perform Legendrian surgeries on the successive Legendrian push-offs of $\alpha_j^i$'s, in the order corresponding to the order of Dehn twists in the decomposition of the monodromy.  

\begin{figure}[htpb!] \begin{center}
\begin{picture}(340,150)
\put(0,0){\includegraphics{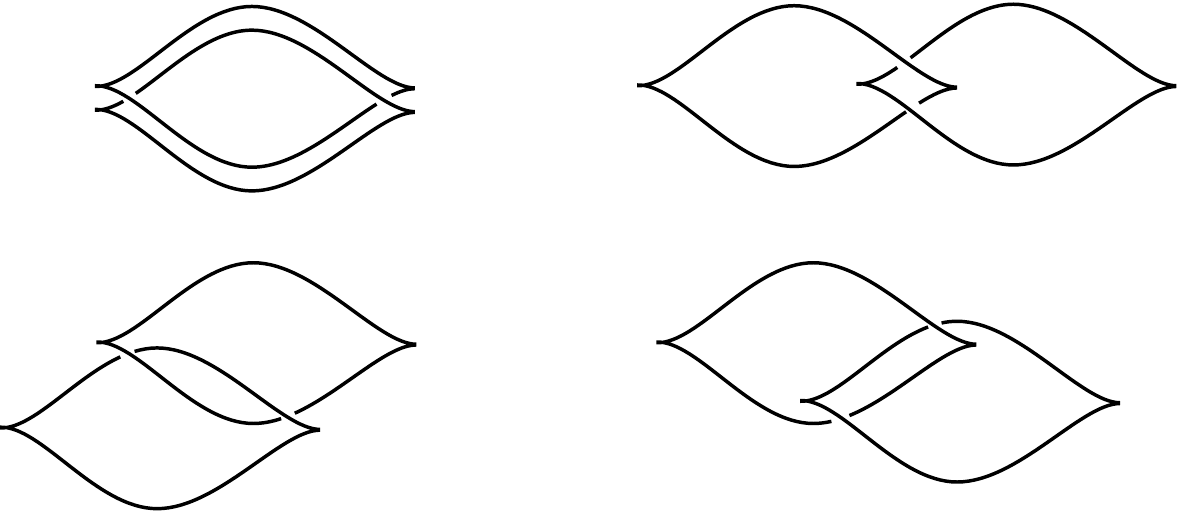}}
\put(110,130){${\alpha_k^j}^+$}
\put(110,100){$\alpha_k^j$}
\put(110, 60){$\alpha_k^j$}
\put(85, 5){${\alpha_{k-1}^{j-1}}^+$}
\put(170, 100){${\alpha_{k-1}^j}^+$}
\put(320, 100){$\alpha_k^j$}
\put(180,30){$\alpha_k^j$}
\put(320,10){${\alpha_k^{j-1}}^+$}
\end{picture}
\caption{Legendrian push-offs of various curves $\alpha^{j}_k$. In all cases not shown, $\alpha^{j}_k$ and  ${\alpha^{i}_l}^{+}$ do not link.}\label{Hopf-links} 
\end{center}\end{figure}
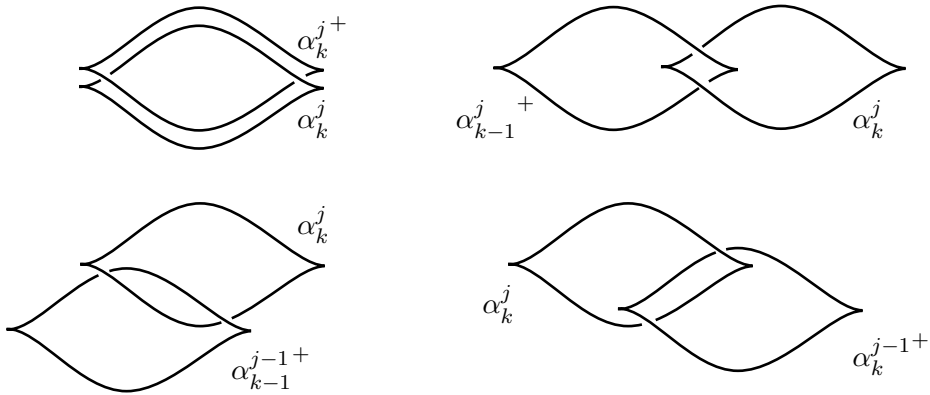

\begin{remark} \label{rmk3.3} In certain cases,  it easy to see that the push-offs of different curves  $\alpha_j^i$ will
be unlinked even if the curves themselves intersect on the surface $P_0$. Indeed, 
consider the braid $K=\sigma_1 \dots \sigma_{n-1}$ and the braid 
$K'=\sigma_1 \dots \sigma_{n-1} \sigma_j$ which differs from $K$ by an additional crossing.
The links $K$ and $K'$ differ only in a small ball $B$ that contains this crossing; 
the $p$-fold branched  cover of $B$ is a genus $(p-1)$ handlebody, and the contact manifolds $\Sigma_p(K)=S^3$ and  $\Sigma_p(K')$ differ only by a surgery on this handlebody. (In fact, the surgery on the handlebody is equivalent 
to $(p-1)$ surgeries on the push-offs of $\alpha_k^j$ where $k=1, \cdots, p-1$, corresponding to the given crossing; 
the surgery curves are all contained in the handlebody.) We also observe that $B$ can be thought of as neighborhood
of the arc connecting two strands of the $K$ at the given crossing, and that the  $p$-fold branched  cover of $B$
is then a neighborhood of the branched $p$-fold cover of this arc. Now, let $c_1$ and $c_2$ be two extra 
crossings added to $K$, and $a_1$ and $a_2$ the corresponding arcs. Untwisting the unknot $K$, we can easily determine 
whether the lifts of $a_1$ and $a_2$ to the branched cover are linked; if they are not, the corresponding surgery
curves will not be linked either.  If, however,  $a_1$ and $a_2$ are linked, we have to examine the push-offs 
of the related curves $\alpha_k^j$ to determine the surgery link. 

We orient $\alpha^{j}_k$  so that 
it goes from $x_{j}$ to $x_{j+1}$ on the $k$-th
sheet of Figure~\ref{Dehn-twist}. Observe that 
\begin{align*}
lk( \alpha_k^j, {\alpha_k^j}^+) &= lk(\alpha_k^j, {\alpha_{k-1}^{j-1}}^+)=-1, \\ 
lk( \alpha_k^j, {\alpha_{k-1}^j}^+) &= lk(\alpha_k^j, {\alpha_k^{j-1}}^+)=+1.
\end{align*}
(See Figure~\ref{Hopf-links}.) In all other cases, the curves  $\alpha_k^j$ and  ${\alpha^i_l}^+$ do not link each other. 
\end{remark}

Given a transverse $n$-braid $L$, we can always write in the form  $L=(\sigma_1 \sigma_2 \dots \sigma_{n-1}) b$, 
(possibly after multiplying by the trivial word $\sigma_1 \dots \sigma_{n-1} \sigma_{n-1}^{-1} \dots \sigma_1^{-1}$).

We will use notation $\Omega_p(L)$ for the contact surgery diagram for the $p$-fold cover $\Sigma_p(L)$ for a transverse $n$-braid $L$, constructed by the above method. Thus $\Omega_p(L)$ is a collection of Legendrian unknots equipped with contact surgery coefficients. When $p$ is fixed, we often drop it from notation. Let $P_\theta$, $\theta \in [0, 2\pi)$ denote the pages. 

Examining the addition of an individual $\sigma_k$ or $\sigma_k^{-1}$ to the braid word for $L$, we obtain the following theorem.  

%%%%%%%%%%%%%%%% Theorem {SurgThm} old & correct %%%%%%%%%%%%%%%%%%%%%%%%%%%

\begin{theorem} \label{SurgThm} Fix $p\geq 2$. Suppose that $L=(\sigma_1 \sigma_2 \dots \sigma_{n-1}) b$
is a transverse $n$-braid, and assume that 
$\Omega_p(L) \subset \displaystyle{\bigcup_{0<\theta < \theta_0}} P_\theta$ for some $\theta_0 < 2 \pi$. 
Pick $\theta_0< \theta_1 < \theta_2 < \cdots < \theta_{p-1} < 2 \pi$. Denote the copy of $\alpha^{i}_j$ in the page $P_\theta$ of the open book by $\alpha^{i,\theta}_j$. 

 (1) Suppose $L^+ = (\sigma_1 \sigma_2 \dots \sigma_{n-1}) b  \sigma_k$,  $1 \leq k \leq n$.
 Define diagram $u^+_k$ as in Figure~\ref{u+}. Then $\Omega(L^+) = \Omega(L) \cup u^+_k$.
\begin{figure}[htpb]
\begin{center}
%\psfrag{+}{$\se{-1}$} 
%\psfrag{2}{$\alpha_1^{k, \theta_{p-1}}$} 
%\psfrag{3}{$\alpha_2^{k, \theta_{p-2}}$} 
%\psfrag{p}{$\alpha_{p-1}^{k,\theta_1}$}  
%\psfrag{q}{$\alpha_{p-2}^{k,\theta_2}$} 
%\includegraphics[height=2cm]{h}
\begin{picture}(238,60)
\put(0,0){\includegraphics{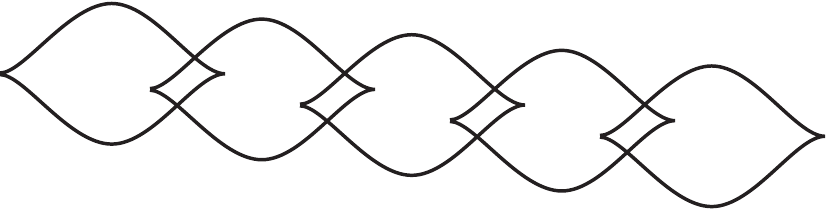}}
\put(30,63){$\se{-1}$}
\put(73,58){$\se{-1}$}
\put(120,52){$\se{-1}$}
\put(164,48){$\se{-1}$}
\put(206,43){$\se{-1}$}
\put(-15,18){$\alpha_1^{k, \theta_{p-1}}$} 
\put(33,9){$\alpha_2^{k, \theta_{p-2}}$} 
%\put(115,52){$\alpha_1^{k, \theta_{p-1}}$}
\put(123,0){$\alpha_{p-2}^{k, \theta_{2}}$}
\put(169,-5){$\alpha_{p-1}^{k, \theta_{1}}$}
\end{picture}
\caption{Contact surgery diagram $u^+_k$.}\label{u+}
\end{center}
\end{figure}

(2) Suppose $L^- = (\sigma_1 \sigma_2 \dots \sigma_{n-1}) b  \sigma_k^{-1}$,  $1 \leq k \leq n$.  
Define diagram $u^-_k$ as in Figure~\ref{u-}. Then $\Omega(L^-) = \Omega(L) \cup u^-_k$. 
\begin{figure}[htpb!]
\begin{center}
\begin{picture}(230,85)
%\psfrag{+}{$\se{+1}$} 
%\psfrag{2}{$\alpha_1^{k, \theta_1}$} 
%\psfrag{3}{$\alpha_2^{k, \theta_2}$} 
%\psfrag{p}{$\alpha_{p-1}^{k,\theta_{p-1}}$}  
%\psfrag{q}{$\alpha_{p-2}^{k,\theta_{p-2}}$} 
\put(0,0){\includegraphics{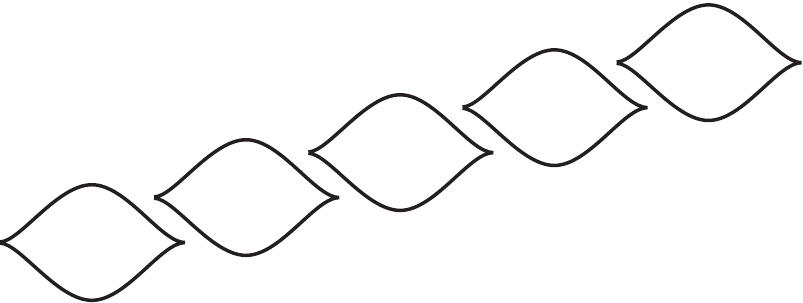}}
\put(-5,0){$\alpha_1^{k, \theta_1}$} 
\put(72,6){$\alpha_2^{k, \theta_2}$} 
\put(168,33){$\alpha_{p-2}^{k,\theta_{p-2}}$}  
\put(216,46){$\alpha_{p-1}^{k,\theta_{p-1}}$} 
\put(24,38){$\se{+1}$}
\put(68,52){$\se{+1}$}
\put(112,66){$\se{+1}$}
\put(158,76){$\se{+1}$}
\put(203,91){$\se{+1}$}
\end{picture}
\caption{Contact surgery diagram $u^-_k$.}\label{u-}
\end{center}
\end{figure}

(Here and below, we draw Legendrian links as their front projections to the $(x, z)$ plane.) 

The diagrams  $u^+_k$ and  $u^-_k$ may be linked  to $\Omega(L)$; the way they link can be determined by drawing 
the corresponding Legendrian push-offs of $\alpha_j^i$ as explained in Remark \ref{rmk3.3}.
\end{theorem}

\begin{proof} This is a direct application of the algorithm developed above.
\end{proof}

\begin{cor} 
Suppose that $L = \sigma_1 \dots \sigma_{n-1} b \in B_{n}$ and 
$L' =  \sigma_1 \dots \sigma_{n-1} \sigma_{n} b \in B_{n+1}$ (i.e., $L'$ is a positive stabilization of $L$ representing the same transverse link).
Then $\Omega(L') = \Omega(L)$.  Note that every positively stabilized braid can be written in such form. 
\end{cor}
\begin{proof} 
The braids $L$ and $L'$ give rise to different initial open books for $S^3$, the one corresponding to the braid 
$\sigma_1 \dots \sigma_{n-1}$ and the other to $\sigma_1 \dots \sigma_{n-1} \sigma_{n}$, but the subsequent 
Dehn twists corresponding to $b$ produce identical surgery diagrams. 
\end{proof}

\begin{cor} \label{cor-ot}
Let $L = \sigma_1 \dots \sigma_{n-1} b \in B_{n}$, and 
$L_{stab} =  \sigma_1 \dots \sigma_{n-1} \sigma_{n}^{-1} b \in B_{n+1}$ 
(i.e., $L_{stab}$ is a {\em negative} braid stabilization of $L$,
representing a transverse link stabilization). Define  $u^{ot}_n$ as in Figure~\ref{u-ot}.
Then the contact manifold represented by $u^{ot}_n$ is an overtwisted 3-sphere, and
$\Omega(L_{stab}):= \Omega(L) \sqcup u^{ot}_n$. 
We use the symbol ``$\sqcup$'' to emphasize that $u^{ot}_n$ is {\em not} linked to $\Omega(L)$ .
\end{cor}
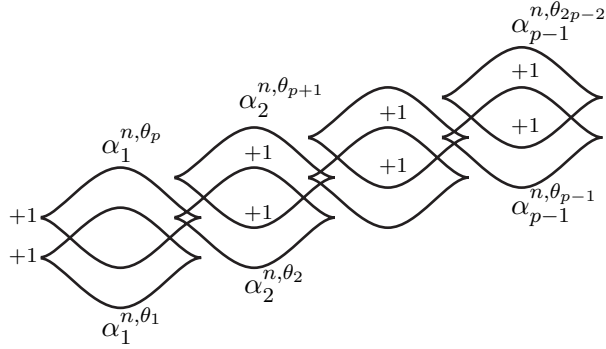
\begin{figure}[htpb!]\label{u-ot}
\begin{center}
%\psfrag{+}{$\se{+1}$} 
%\psfrag{1}{$\alpha_1^{n, \theta_1}$} 
%\psfrag{2}{$\alpha_2^{n, \theta_2}$} 
%\psfrag{4}{$\alpha_{p-1}^{n,\theta_{p-1}}$} 
%\psfrag{5}{$\alpha_1^{n,\theta_p}$} 
%\psfrag{6}{$\alpha_2^{n,\theta_{p+1}}$} 
%\psfrag{8}{$\alpha_{p-1}^{n, \theta_{2p-2}}$}
\begin{picture}(212,110)
\put(0,10){\includegraphics{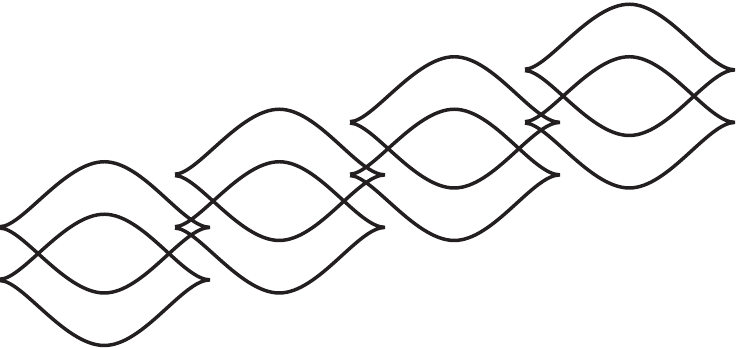}}
\put(23,0){$\alpha_1^{n, \theta_1}$} 
\put(23,69){$\alpha_1^{n,\theta_p}$} 
\put(77,16){$\alpha_2^{n, \theta_2}$} 
\put(75,86){$\alpha_2^{n,\theta_{p+1}}$}
\put(178,46){$\alpha_{p-1}^{n,\theta_{p-1}}$} 
\put(178,115){$\alpha_{p-1}^{n, \theta_{2p-2}}$}
\put(-12,43){$\se{+1}$}
\put(-12,28){$\se{+1}$}
\put(77,44){$\se{+1}$}
\put(77,67){$\se{+1}$}
\put(128,82){$\se{+1}$}
\put(128,62){$\se{+1}$}
\put(178,77){$\se{+1}$}
\put(178,98){$\se{+1}$}
\end{picture}
\caption{Contact surgery diagram $u^{ot}_n$.}
\end{center}
\end{figure}

\begin{proof} We write $L_{stab} =  \sigma_1 \dots \sigma_{n-1} \sigma_{n} \sigma_{n}^{-2} b$ and apply part 2 of Theorem~\ref{SurgThm} twice. To show that the contact manifold represented by 
$u^{ot}_n$ is an overtwisted 3-sphere, we first use  Kirby calculus to see  
that the underlying smooth manifold is $S^3$. Using formula (\ref{d3}), 
we compute $d_3=-\frac12 + p-1$. (We have $c_1(s_J)=0$, $\sign(X)=0$, $\chi(X)=1+2(p-1)$.) 
Since we know that $\xi_{std}$ is the unique tight contact structure on $S^3$, and $d_3(\xi_{std})=-\frac12$, 
it follows that the contact structure given by  the diagram $u^{ot}_n$
is overtwisted. The branched cover of $L_{stab}$ is then the connected sum of this overtwisted sphere 
and the branched cover of $L$.
\end{proof}

\begin{cor}
Suppose that $L_n = \sigma_1 \dots \sigma_{n-1} b \in B_{n}$ is an $n$-braid, 
and $L_{n+1} = \sigma_1 \dots \sigma_{n-1} b \in B_{n+1}$
is an $(n+1)$-braid obtained from $L_n$ by an addition
of a trivial $(n+1)$-th strand. Then $\Omega(L_{n+1}) = \Omega(L_n) \sqcup u^-_n$, where  
$u^-_n$ is {\em not} linked to $\Omega(L)$.
\end{cor}

\begin{proof} This follows from the identity $L_{n+1} = \sigma_1 \dots \sigma_{n-1} \sigma_n \sigma_n^{-1} b$, 
and the fact that the word $b$ does not contain $\sigma_n^{\pm n}$. We also observe that on the level of contact
manifolds, we are taking a connected sum with $\#_p S^1\times S^2$, where the latter is equipped with its unique 
Stein fillable contact structure.   
\end{proof}

\begin{figure}[htpb!]
\begin{center} 
\includegraphics[scale=0.6]{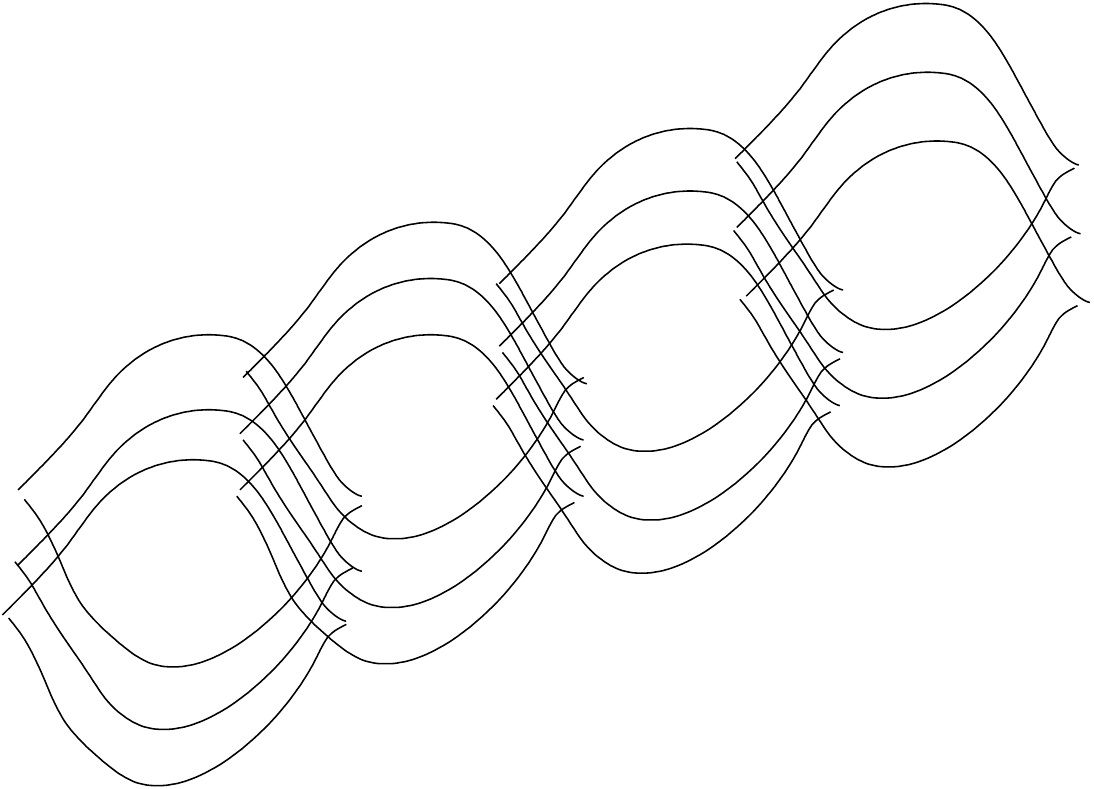}
\caption{A surgery diagram for $\Sigma_5(\sigma_1^4)$. We perform Legendrian surgery on each component.}\label{surgery-example} 
\end{center}
\end{figure}

It is now easy to obtain surgery diagrams of all $p$-fold branched covers of 2-braids. 
\begin{example} 
A surgery diagram for the 5-fold cover of the transverse braid $\sigma_1^4$ is shown in Figure~\ref{surgery-example}. 
\end{example}

\begin{remark}\label{rem-different-choice}

Even though every closed $n$-braid is isotopic to a braid containing 
a string $\sigma_1 \sigma_2 \sigma _3 \dots \sigma_{n-1}$, we may want to start with 
an open book corresponding to another version of transverse unknot, 
say $\sigma_2 \sigma_1 \sigma_3 \dots, \sigma_{n-1}$ (this will be useful in Subsection~\ref{subsec-not}).
To obtain this other open book, we consider Figure~\ref{LegReal} and change the curves 
$\alpha_1^2, \alpha_2^2, \dots, \alpha_{p-1}^2$, so that they now go through the top two rows 
of the grid-like page. This is shown on Figure~\ref{AnotherLegReal}; the other curves 
$\alpha_1^k, \alpha_2^k, \dots \alpha_{p-1}^k$ for $k\neq 2$ remain the same.
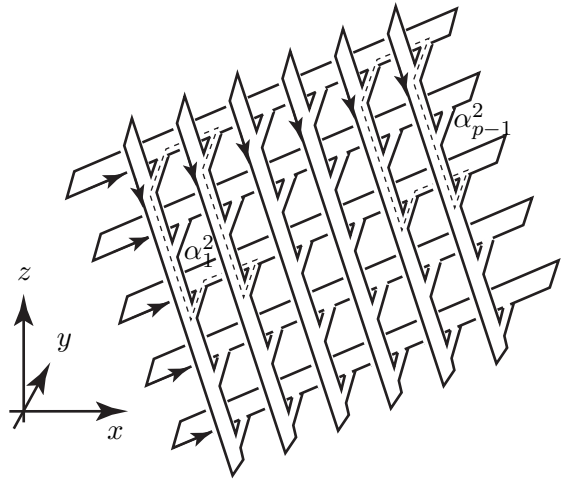
\begin{figure}[htpb!] \begin{center}
\begin{picture}(210, 185)
\put(0,0){\includegraphics{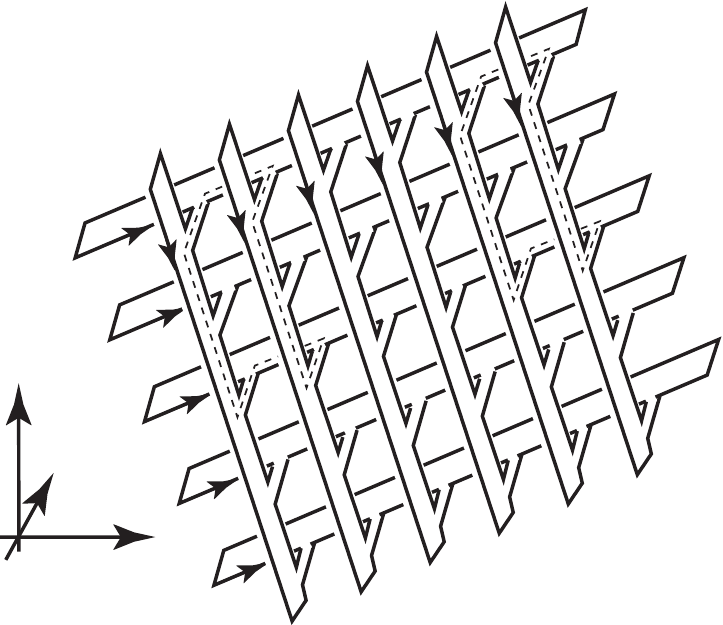}}
\put(37, 15){$x$}
\put(18, 50){$y$}
\put(3,75){$z$}
\put(66, 83){$\alpha^2_1$}
\put(168, 133){$\alpha^2_{p-1}$}
\end{picture}
\caption{A different choice of curves $\alpha_j^2$ produces an open book 
whose monodromy is $\sigma_2 \sigma_1 \sigma_3 \dots \sigma_{n-1}$.}\label{AnotherLegReal} 
\end{center}
\end{figure}  
For open books, this change corresponds to plumbing positive Hopf bands together 
in a slightly different way to form the same page. Analyzing the push-offs of the curves $\alpha_j^k$ as 
in Proposition \ref{sigma1}, we see that the monodromy of the open book can now be expressed 
as  $(D^{n-1}_1 \circ \cdots \circ D^{n-1}_{p-1}) \circ \cdots \circ  (D^{1}_1 \circ \cdots \circ D^1_{p-1})\circ(D^2_1 \circ \cdots \circ D^2_{p-1})$, which by Lemma~\ref{Dehn twist lemma} corresponds to the braid 
$\sigma_2 \sigma_1 \sigma_3 \dots \sigma_{n-1}$ as required.

Another case worth mentioning is the initial unknot given by the braid $\sigma_n \dots \sigma_2 \sigma_1$.
In this case, we have the same open book as for the unknot $\sigma_1 \sigma_2  \dots \sigma_n$,
with the role of the curve $\alpha_j^{k}$ played by  $\alpha_j^{n-k}$.

In principle, it is not necessary to single out the braid word that gives the unknot: we can as well
start from the trivial braid and obtain $(\Sigma(L),\xi_L)$ as a result of surgery
on ${\#_{p-1}} S^1\times S^2$. However, the presence of 1-handles seems to complicate 
matters.
\end{remark}

%%%%%%%%%%%%%%%%%%%%%%%%%%%%%%%
\section{Properties of Branched Covers} \label{properties}

In this section, we  prove Theorems~\ref{tight-ow} and \ref{homotp}.
The proofs are very similar to those of \cite[sections 4 and 5]{Pla}.

%-----------------------------------------------------------------
\subsection{Quasipositive braids and Stabilizations.}

Recall \cite{Ru} that a braid is called {\em quasipositive} if its braid
word is a product of conjugates of the standard generators.

 \begin{proof}[Proof of Theorem~\ref{tight-ow}]
 If $L$ is quasipositive, we can resolve a few positive crossings to convert the braid representing $L$ into a braid
equivalent to a trivial one (of the same braid index). The $p$-fold cover branched 
over the trivial braid is a connected sum of several copies of
$(S^1\times S^2, \xi_0)$, which is Stein fillable ($\xi_0$ here stands for the unique Stein fillable contact structure 
on $S^1\times S^2$).
Putting  the positive crossings back in,  
by Lemma~\ref{Dehn twist lemma}  we see that the monodromy of the open book
for $(\Sigma_p(L), \xi_p(L))$ is given by a composition of positive 
Dehn twists. It follows that the contact manifold is Stein fillable. 
 The second part of the theorem follows from Corollary~\ref{cor-ot}.
\end{proof}

%-----------------------------------------------------------------
\subsection{Homotopy Invariants.}

\begin{proof}[Proof of Theorem~\ref{homotp}]
The fact that $c_1(\s_{L})=0$ follows immediately:  $\s_\xi$ is the restriction to $Y$ of
the $\Spinc$ structure  $\s_J$ described in Subsection \ref{cont-surgery}; $c_1(\s_J)$ evaluates as
$0$ on each homology generator corresponding to either a $(-1)$ or
a $(+1)$ surgery, because  all surgeries are performed 
on standard Legendrian unknots with rotation number $0$.

For the second part of the theorem, 
suppose that two closed braids $L$ and $L'$ are
isotopic as smooth knots, and that $sl(L)=sl(L')$.  By the Markov theorem for smooth knots 
\cite{Bi}, $L'$  can be obtained from $L$ by a sequence of braid isotopies and (positive
and negative) braid stabilizations and destabilizations.
Braid isotopies and positive stabilizations preserve both $sl$ 
and the $d_3$ invariant, since they do not change the transverse link type.
Each negative stabilization (resp. destabilization) decreases (resp. increases)
the self-linking number by $2$ and the $d_3$ invariant by $p-1$, 
since, as we saw in Corollary~\ref{cor-ot}, transverse stabilization gives the connected sum with the overtwisted sphere in Figure~\ref{u-ot}. But if  $sl(L)=sl(L')$, every negative stabilization must be compensated by a negative destabilization. It follows that $d_3(\xi_L)= d_3(\xi_{L'})$.  
\end{proof}

\begin{cor} 
Fix $p\geq 2$. Let $T$ be  a transverse link smoothly isotopic to
the $(m, n)$ torus link, $m, n>0$. The branched cover $\Sigma_p(T)$ is then the Brieskorn 
sphere $\Sigma(m, n, p)$.
 If $sl(T)=sl_{max}=mn-m-n$, then
$\xi_p (T)$ is Stein fillable. Otherwise 
$\xi_p (T)$ is overtwisted. For different values of $sl(T)$, these
overtwisted structures have different $d_3$ invariants. 
\end{cor}

\begin{proof} 
By classification of transverse torus links, $T$ is given by a positive 
braid when it has the maximal self-linking number; otherwise it is isotopic
to this link transversely stabilized $r$ times (for some $r>0$).
The  number of stabilizations $r$ can be determined by the self-linking number. 
\end{proof}

%-----------------------------------------------------------------
\subsection{Overtwisted branched covers.}

We generalize the second part of Theorem~\ref{tight-ow}, and show that the branched covers are overtwisted for a large family of transverse links.

\begin{prop} \label{neg-crossings} 

Suppose that the transverse link $L$ is represented by a
closed braid such that its braid word $\phi_L$ contains a factor of
$\sigma_i^{-1}$ but no $\sigma_i$'s for some $i>0$. (This means
that all the crossings in the braid diagram on the level between
the $i$-th and the $(i+1)$-th strands are negative.)
Then the branched 
$p$-fold cover $\Sigma_p(L)$ is overtwisted for any $p\geq 2$.

\end{prop}

\begin{proof} We will use the right-veering monodromy criterion for tightness \cite{HKM}.

More precisely, we will show that the monodromy of the open book for the branched cover of $L$ given by Lemma \ref{Dehn twist lemma} is not right-veering, which by \cite[Theorem 1.1]{HKM} implies that the corresponding 
contact structure is overtwisted. (This is in fact the criterion for overtwistedness that was first 
given in \cite{Goo} in terms of ``sobering arcs''.) 

The negative crossings between $i$-th and $(i+1)$-th strands yield the open book $({\tilde D}, {\tilde \phi}_L)$ such that the Dehn twists about curves $\alpha^{i}_j$ are all left-handed for $j=1, \dots, p-1$. 

\begin{figure}[htpb!] 
\begin{center}
\includegraphics[scale=0.8]{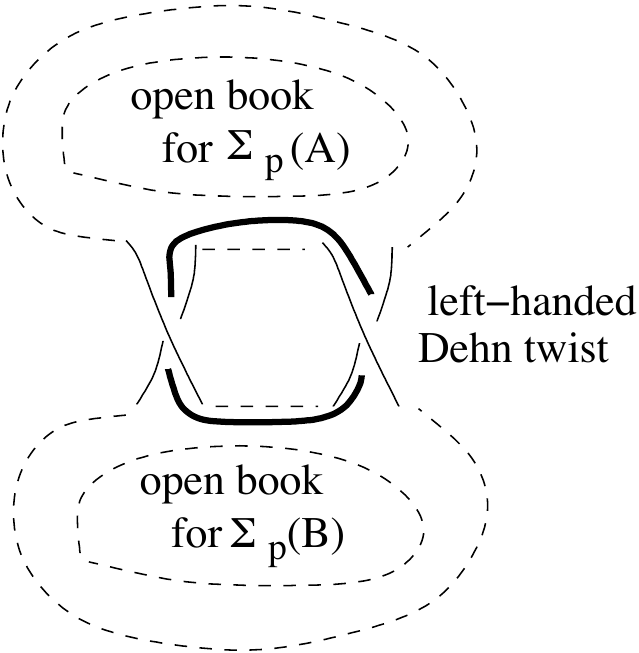}
\caption{The plumbed sum of the open books $({\tilde D}, {\tilde \phi}_A)$, $({\tilde D}, {\tilde \phi}_B)$, and 
a left-handed Hopf band.} \label{destab}
\end{center}
\end{figure}  

Remove from $L$ all the (negative) crossings between  the $i$-th and the $(i+1)$-th strands (in other words, 
remove all the negative factors of $\sigma_i$ from the braid word  $\phi_L$). The link $L$ then splits into two 
links $A$ and $B$ such that the corresponding braid words $\phi_A$ and $\phi_B$ contain only generators $\sigma_j$
with $j<i$ resp. $j>i$. Let    $({\tilde D}, {\tilde \phi}_A)$,   $({\tilde D}, {\tilde \phi}_B)$ be the 
open books for the branched covers $\Sigma_p(A)$ and $\Sigma_p(B)$, and 
consider the plumbed sum of  $({\tilde D}, {\tilde \phi}_A)$, $({\tilde D}, {\tilde \phi}_B)$, and 
a left-handed Hopf band. The resulting open book is shown on Figure~\ref{destab}, and is in fact 
a negative stabilization of the connected sum $({\tilde D}, {\tilde \phi}_A)\#({\tilde D}, {\tilde \phi}_B)$.
It is easy to see that the monodromy of a negatively-stabilized open book is not right-veering. 
The open book $({\tilde D}, {\tilde \phi}_L)$ is obtained from this non-right-veering open book by a sequence 
of negative stabilizations and additional left-handed Dehn-twists, and so cannot be right-veering either
(because a composition of a right-handed Dehn twist and a right-veering monodromy is right-veering). 
\end{proof}

\begin{remark}
If the branched $p$-fold covers of two transverse links $L_1$ and $L_2$ of the same topological type are both overtwisted and $sl(L_1)= sl(L_2)$, then Theorem~\ref{homotp} together with Eliashberg's classification of overtwisted contact structures implies that $\Sigma_p(L_1)$ is contactomorphic to $\Sigma_p(L_2)$. 
\end{remark}

We therefore have

\begin{cor} \label{bmcont} 
In Table 1 in \cite{BM1} of transverse knots, all pairs (except perhaps for the representatives of the knot $11_{a240}$) give rise to contactomorphic $p$-fold branched covers for all $p\geq 2$.  
\end{cor}

In view of the previous remark, showing that certain branched covers are overtwisted can be useful.
We thus illustrate two other ways to establish overtwistedness (our examples below are all included 
in Proposition \ref{neg-crossings}, but the methods can be used for other links as well).

The first method applies in the rare cases where the classification of tight contact structures is known 
for the smooth manifold $\Sigma_p(L)$. For example, this is the case for double covers of  2-bridge links:
it is well known that these are lens spaces, and the tight contact structures on lens spaces were classified 
in \cite{Ho}. 

\begin{figure}[ht] 
\begin{center}
\includegraphics[scale=0.8]{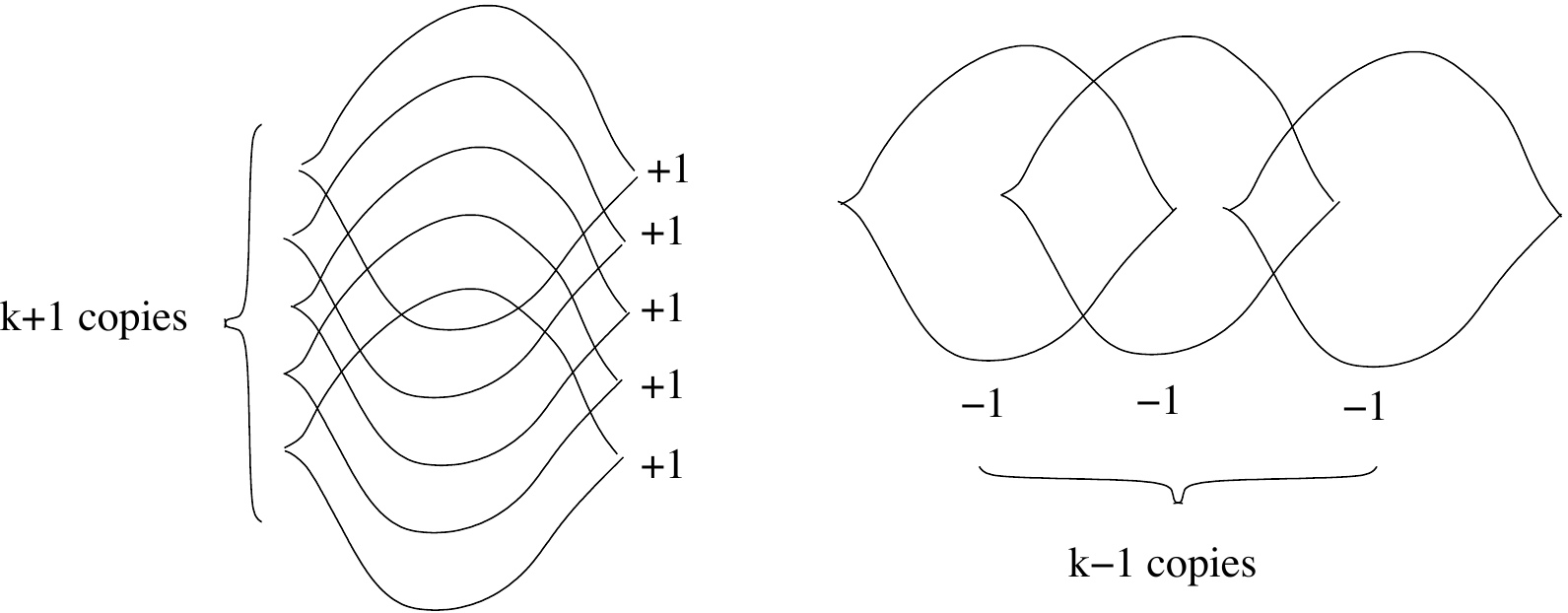}
\caption{The branched double cover of $\sigma_j^{-k}$ (left) and the unique tight contact structure on $L(k, k-1)$ (right).} \label{k-neg}
\end{center}
\end{figure}

Consider the transverse 2-braid $L=\sigma_1^{-k}$ where $k \geq 1$; its branched double cover is the lens space 
$-L(k,1)= L(k, k-1)$,  with the contact structure $\xi$ given by the surgery diagram on Figure~\ref{k-neg}
(where $(+1)$ contact surgery  is performed  on each of $k+1$ successive push-offs of the Legendrian unknot of $tb=-1$). We compute the $d_3$ invariant of this contact structure. If $X$ is the 4-manifold corresponding to the surgery, then $\sign(X)=k-1$ (indeed, the intersection form for $X$ has zeroes on the diagonal and $-1$'s for all other entries; it is easy to see that the matrix has an eigenvalue $1$ of order $k$ and an eigenvalue $-k$ of order 1). We also have $c_1(X)=0$, and $\chi(X)=k+2$. Therefore, from (\ref{d3}) we obtain
$d_3(\xi)= \frac{-k+3}4$. On the other hand, by \cite{Ho}, the lens space $L(k, k-1)$ carries a 
unique tight contact structure $\xi_0$; this contact structure is the boundary of a linear plumbing (also shown on Figure~\ref{k-neg}). The corresponding Stein 4-manifold $X_0$ has $c_1(X_0)=0$, $\sign(X_0)=0$ and $\chi(X_0)=k$, so 
$d_3(\xi_0)=-\frac{k}2$. It follows that the contact structure $\xi$ is not isotopic to $\xi_0$, 
and therefore must be overtwisted.

\begin{figure}[htpb] \begin{center}
\begin{picture}(288, 144)
\put(0,0){\includegraphics%[height=4.7cm]
{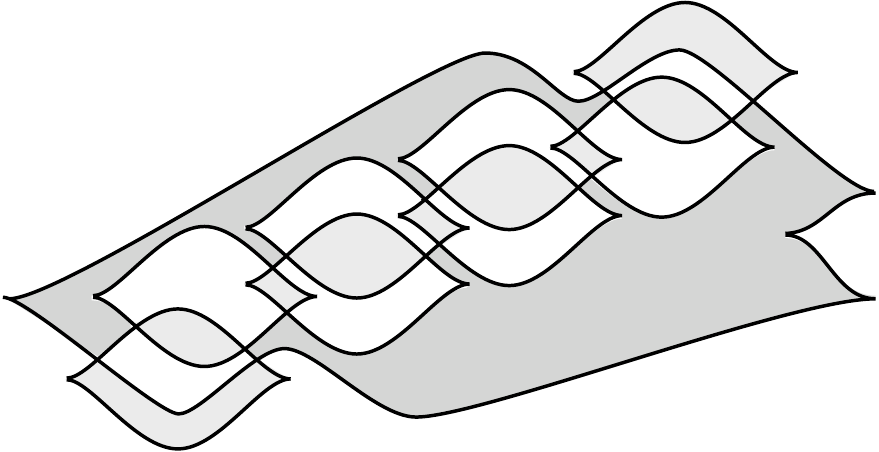}}
\put(90, 100){$K$} \put(213, 46){$S$}
\end{picture}
\caption{The Legendrian unknot knot $K$  bounds an overtwisted disk in the surgered manifold.}\label{otdisk}
\end{center}\end{figure}

Another way to prove overtwistedness is simply to find an overtwisted disk in the surgery diagram. Admitting that these pictures get unwieldy even for simple links, we exhibit such a disk for the overtwisted sphere $u^{ot}$ described in Figure~\ref{u-ot} (i.e. the branched $p$-fold cover of $\sigma^{-1}$). Indeed,  the surface $S$ shown on Figure~\ref{otdisk} induces the $0$-framing on each component of the surgery link $u^{ot}$, and the $(-2)$-framing on the Legendrian knot $K$. (We assume that all Legendrian knots are oriented as the boundary of $S$). Since $(+1)$-contact surgery is $0$-framed Dehn surgery, $S$ becomes a disk bounded by $K$ in the surgered manifold. Then the equality
$tb(K)= \text{ the surface framing of } K = -2$ implies that this disk is overtwisted.

%%%%%%%%%%%%%%%%%%%%%%%%%%%%%%%%%

\section{ Can we distinguish transverse knots?} \label{distinguish}

We can now use the constructions from previous sections to examine the branched covers 
of certain transverse knots and prove Theorem~\ref{introthm} (see Corollaries~\ref{cor-A}, \ref{cor-B}, and Theorem~\ref{neg-flype-thm}). We already saw that for most examples of \cite{BM1}, 
the branched covers do not detect the difference between transverse knots.
We now consider the remaining pairs of non-isotopic transverse
knots with the same classical invariants from \cite{BM}, \cite{BM1}, \cite{NgOT}, 
and try to distinguish them via the corresponding contact structures. 

\subsection{Birman--Menasco examples}

The methods of Birman-Menasco \cite{BM}, \cite{BM1} produce examples that are pairs of 3-braids $L_1, L_2$ related by a {\em negative flype}. This means that $L_1 = \sigma_1^u \sigma_2^v  \sigma_1 ^w  \sigma_2^{-1}$ and  $L_2 = \sigma_1^u  \sigma_2^{-1} \sigma_1 ^w \sigma_2^v$.  

Recall that the contact structure $\bar \xi$  {\em conjugate  to} $\xi$ is obtained from $\xi$ by reversing the orientation of contact planes.  
 
\begin{prop} Transverse 3-braids $L_1$ and $L_2$ related by a negative flype give rise to conjugate contact structures on the branched covers: $\xi_p(L_1)$ is contactomorphic to $\bar{\xi}_p(L_2)$. \label{3braids}
\end{prop}

\begin{figure}[htb!]
\begin{center} 
\includegraphics[scale=1.0]{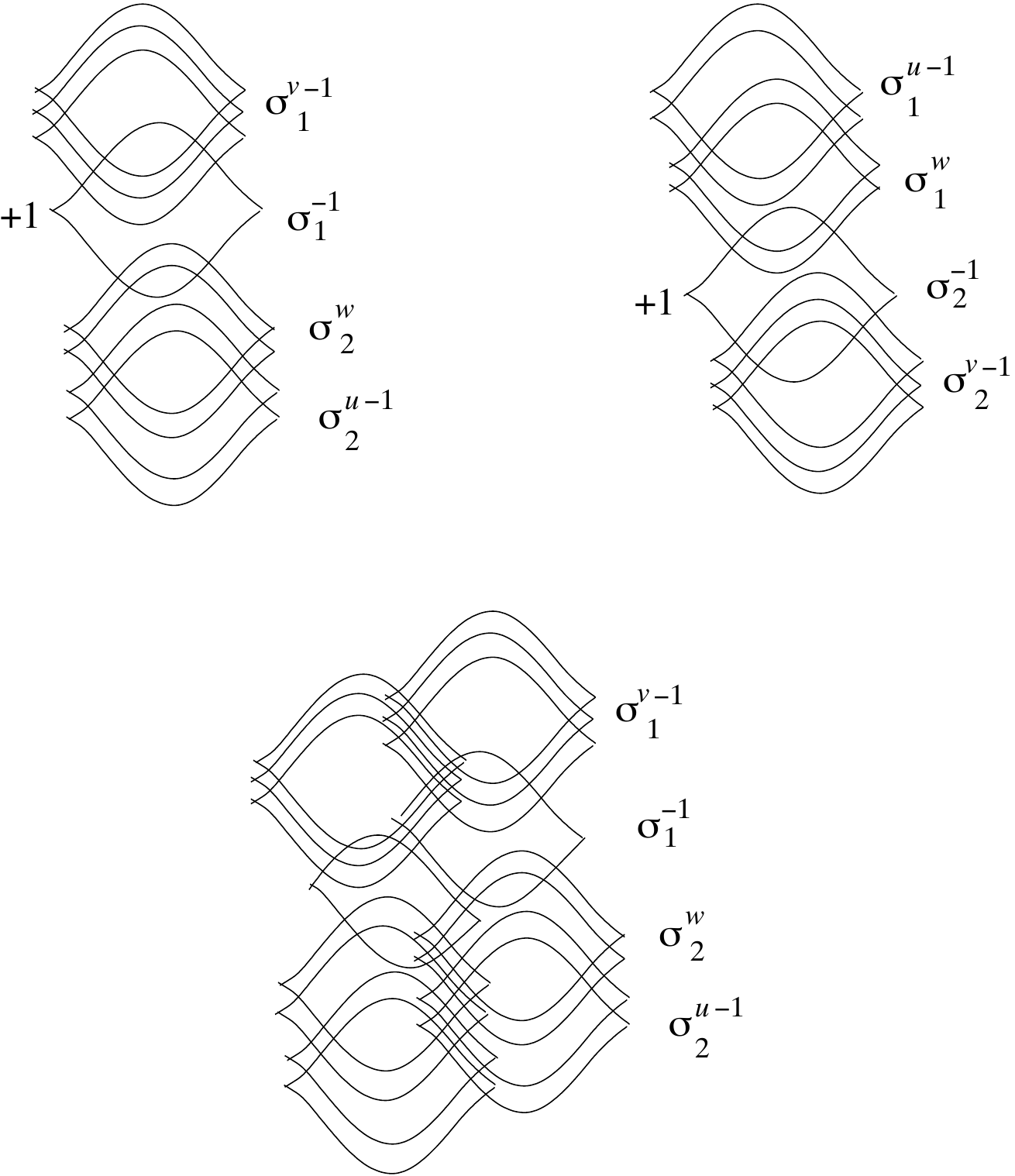}
\caption{Branched double covers of  $L_1$ (top left) and $L_2$ (top right). We assume that $u-1, v-1, w \geq 1$;
a $(+1)$ contact surgery is performed on unknots marked with $+1$, Legendrian surgery on all other components.
A diagram for the branched 3-fold cover of $L_2$ is shown at the bottom; branched 3-fold cover of $L_1$ is obtained by contact surgery on a mirror of this link. To obtain $p$-fold covers, take $(p-1)$ copies of the surgery link 
for the double cover linked in a way similar to the 3-fold cover case.}    \label{L1L2cover}
\end{center}
\end{figure}

\begin{proof} 
We write the closed braids as 
$$
L_1=(\sigma_1  \sigma_2) \sigma_2^{v-1} \sigma_1^w \sigma_2^{-1} \sigma_1^{u-1}, \quad L_2 = (\sigma_2 \sigma_1) \sigma_1^{u-1} \sigma_2^{-1} \sigma_1^w \sigma_2^{v-1}.
$$
Observe that $L_2$ can be taken to  
$(\sigma_1 \sigma_2) \sigma_2^{u-1} \sigma_1^{-1} \sigma_2^w \sigma_1^{v-1}$ by a transverse isotopy.
Using the method in Theorem~\ref{SurgThm} and the following corollaries, we can draw surgery diagrams for the branched covers of $L_1$ and $L_2$.
For example, double covers for the case where $u-1, v-1, w \geq 1$ are shown on Figure~\ref{L1L2cover} (top);
we see that they are obtained by contact surgeries on two links which are Legendrian mirrors of one another.
Similarly, $p$-fold branched covers for $L_1$ and $L_2$ are also obtained by surgery on Legendrian mirrors,
since the corresponding diagrams are obtained by taking $(p-1)$ copies of the surgery link for the double cover
linked as dictated by Figure~\ref{Hopf-links}.
For example, the triple cover for $L_2$ is shown on Figure~\ref{L1L2cover} (bottom).
 For negative $u$, $v$ or $w$ the pictures are similar; besides, the case $v\leq 0$ is already covered by Corollary~\ref{bmcont}. 

We have shown that the surgery link diagram $\Omega_p(L_2)$ for $L_2$ is the Legendrian mirror of the link $\Omega_p(L_1)$ for $L_1$. Now, observe that one link is taken to the other by the map $(x, y, z) \mapsto (-x, y, -z)$. This map reverses the sign of the standard contact form $dz -y\, dx$ (i.e. the orientation of contact planes on $S^3$)
and extends to the map of branched covers that takes $\xi_p(L_1)$ to  $\bar {\xi} _p(L_2)$.   
\end{proof}

\begin{remark} 
Alternatively,  the previous proposition can be proved by using open books. A careful examination of the monodromy shows that the open book for $(\Sigma_p(L_2), \xi_p(L_2))$ can be obtained from the open book for $(\Sigma_p(L_1), \xi_p(L_1))$ by reversing the orientation of the pages as well as the orientation of the $S^1$ direction in the mapping torus. This operation preserves the orientation of the 3-manifold but reverses the orientation of contact planes.
\end{remark}

We can generalize Proposition~\ref{3braids} as follows.

\begin{prop} Let the braid $L_2$ be obtained by reading the braid word $L_1$ backwards, i.e., if
$L_1 = \sigma_{i_1} \sigma_{i_2} \cdots \sigma_{i_{k-1}} \sigma_{i_k}$ then $L_2  =  \sigma_{i_k} \sigma_{i_{k-1}} \cdots \sigma_{i_2} \sigma_{i_1}$.  
Then the contact structures $\xi_p(L_1)$ and  $\xi_p(L_2)$ are conjugate to one another for any $p\geq 2$.
\end{prop}

\begin{proof}
Write 
$$
L_1= (\sigma_1 \sigma_2 \dots \sigma_{n-1})  \sigma_{j_1} \sigma_{j_2} \dots \sigma_{j_{l-1}} \sigma_{j_l}.
$$
Then the braid word for $L_2$ is conjugate to 
$$
L_2=(\sigma_{n-1}  \dots \sigma_2 \sigma_1)  \sigma_{j_l} \sigma_{j_{l-1}} \dots \sigma_{j_2} \sigma_{j_1}. 
$$
In the surgery diagram for cover of $L_1$, the part of the surgery link 
corresponding  to $\sigma_{j_r}$ will be above that for  $\sigma_{j_s}$ when $r<s$;
for cover of $L_2$, it will be below.
In both cases, the surgery unknots corresponding to $\sigma_{j_r}$ and $\sigma_{j_s}$ with $r<s$ will 
be linked (in exactly the same way) iff $j_r \leq j_s$; using the braids-to-surgeries description 
from Section~\ref{surg for cov}, cf. Figure~\ref{Hopf-links}, we see that in fact the surgery links for the two 
branched covers are Legendrian mirrors of one another. It follows that the resulting 
contact structures $\xi_p(L_1)$ and  $\xi_p(L_2)$ are conjugate to one another.
  
(Alternatively, we could rotate  $L_2$ to get 
$$
L_2=(\sigma_1  \dots \sigma_{n-2} \sigma_{n-1})  \sigma_{n-j_l} \sigma_{n-j_{l-1}} \dots \sigma_{n-j_2} \sigma_{n-j_1}, 
$$
and draw the surgery diagrams similar to the Birman-Menasco braids in Proposition~\ref{3braids}.)
\end{proof}

\begin{prop} For any  transverse link $L$, $p\geq 2$, the contact structure ${\xi}_p(L)$ is isomorphic to its conjugate $\bar{\xi}_p(L)$.    
\end{prop}

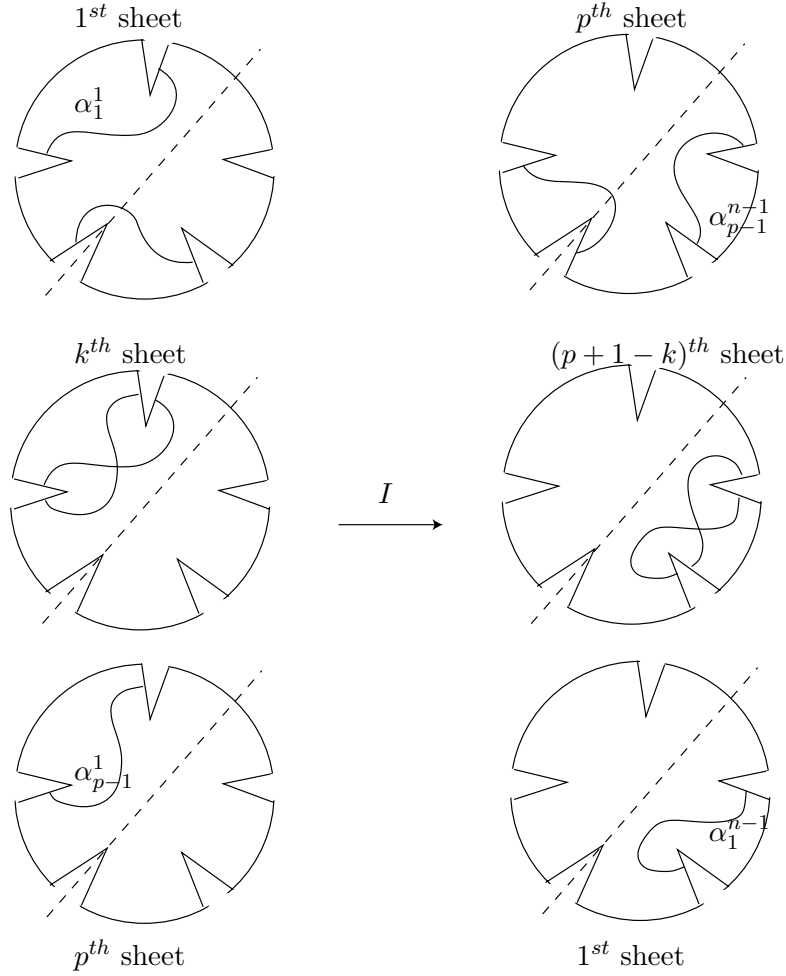
\begin{figure}[htpb] \begin{center}
\begin{picture}(309,354)
\put(0,0){\includegraphics%[height=12.6cm]
{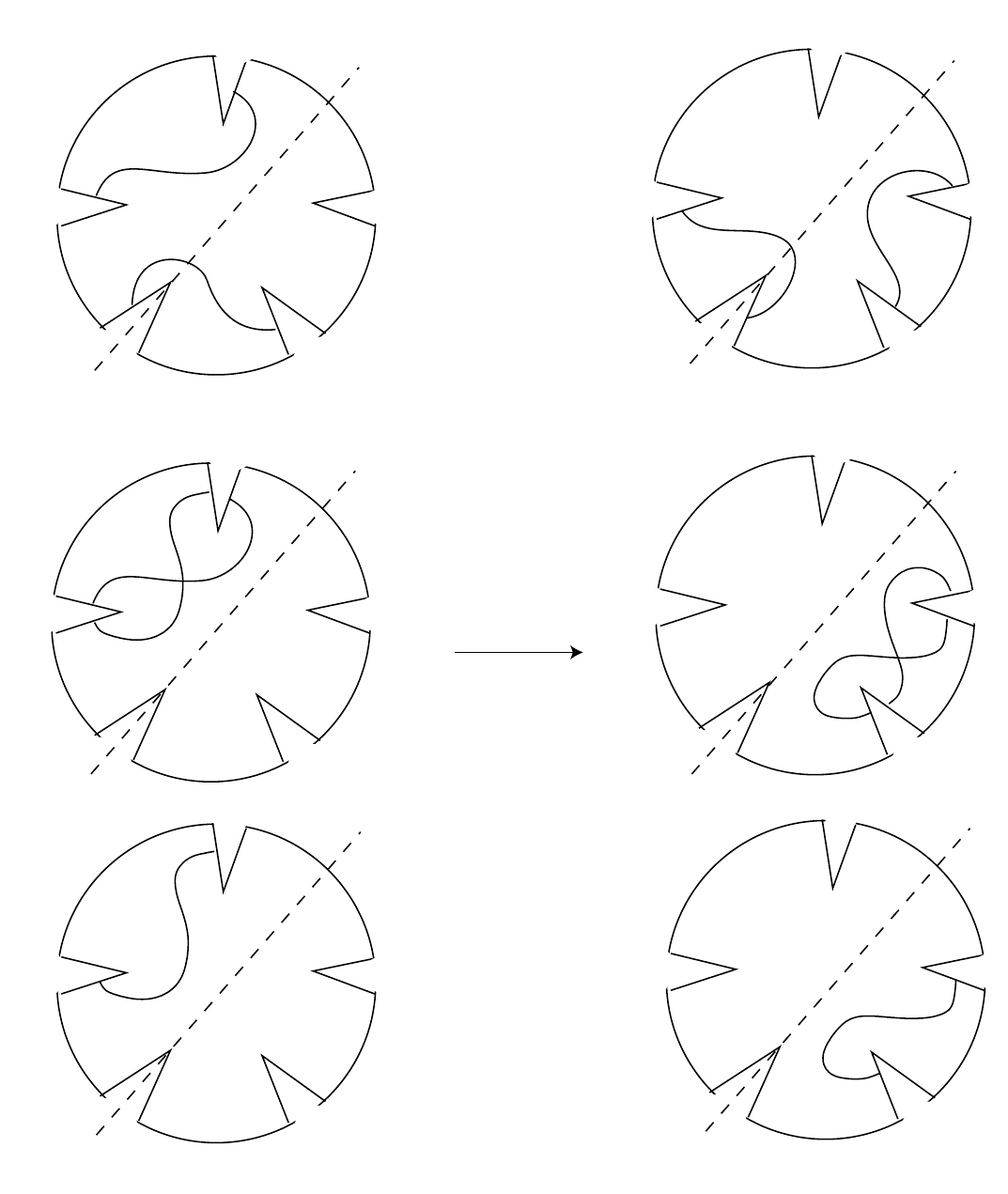}}
\put(40, 355){$1^{st}$ sheet} \put(230, 355){$p^{th}$ sheet}
\put(40,228){$k^{th}$ sheet} \put(220,228){$(p+1-k)^{th}$ sheet}
\put(40, 0){$p^{th}$ sheet}  \put(230,0){$1^{st}$ sheet}
\put(40, 323){$\alpha_1^1$} \put(280, 280){$\alpha_{p-1}^{n-1}$}
\put(40, 70){$\alpha_{p-1}^1$} \put(280, 47){$\alpha_1^{n-1}$}
\put(155, 175){$I$}
\end{picture}
\caption{The involution $I$ on a page of the open book.}\label{reflect} 
\end{center}\end{figure}  

\begin{proof} 
We need to find an involution of the smooth manifold $\Sigma_p(L)$ that induces the orientation 
reversal on contact planes. For a page $P$ of the open book described in Lemma~\ref{Dehn twist lemma}, 
there is an orientation-reversing map $I: P \to P$ that maps $k$-th sheet to the $(p+1-k)$-th sheet, acting as a reflection, and takes the curve $\alpha_k^j$ to the curve $\alpha_{p-k}^{n-j}$ (see Figure \ref{reflect}). If $p$ is odd, the $(p+1)/2$-th sheet is 
mapped to itself, and if $n$ is even, the curve $\alpha_{(p+1)/2}^{n/2}$ is mapped to itself.
Moreover,  $(D_k^j)^{-1} I = I D_{p-k}^{n-j}$, i.e. the involution $I$ takes right-handed Dehn twists 
to the left-handed ones. If $\tilde \sigma_j$ is the lift of the half-twist $\sigma_j$ as in Figure~\ref{sigma}, we have 
$$
(\tilde \sigma_j)^{-1} I = (D_{p-1}^j)^{-1} \dots (D_{2}^j)^{-1} (D_{1}^j)^{-1} I = I D_1^{n-j} D_2^{n-j} \dots D_{p-1}^{n-j} = I \tilde \sigma_{n-j}.
$$ 

Write $\phi_L$ for the braid word for $L$, and let $\phi_{L'}$ be the braid word obtained by changing 
every half-twist generator $\sigma_j$ to $\sigma_{n-j}$. The braids $\phi_L$  and $\phi_{L'}$ are related 
by a braid isotopy (rotating the braid), so if  $L'$ is the transverse link corresponding 
to the braid $\phi_{L'}$, then $L$ and $L'$ are transversely isotopic. However, we have
$$
(\tilde \phi_L)^{-1} I = I \tilde \phi_{L'}. 
$$ 

If we extend the map $I$ to an orientation-preserving  map 
$R: P\times [0,1] \to P\times [0,1]$, defined by $R(x, t) = (I(x), 1-t)$, 
$R$ descends to open books, taking the open book $(P, \phi_L)$ to the open book  $(-P, (\phi_{L'})^{-1})$.
The latter open book is compatible with the contact structure $\bar \xi_{L'}$ which is isotopic to 
$\bar \xi_{L}$. It follows that $\xi_L$ and  $\bar \xi_{L}$ are isomorphic.
\end{proof}

The last two propositions apply in the following special cases, proving Theorem~\ref{introthm}.

\begin{cor}\label{cor-A} 
Let $L$ be a Legendrian link, $\bar L$ its Legendrian mirror, and consider the transverse push-offs $L^+$ and $\bar L^-$. Then the corresponding $p$-fold  branched covers are contactomorphic for all $p$. 
\end{cor}

\begin{cor}\label{cor-B} 
If  3-braids $L_1$ and $L_2$ are related by a negative flype, then $\xi_p(L_1)$ and  ${\xi}_p(L_2)$ are isomorphic.
\end{cor}
 
\begin{remark} The double branched covers of the Birman--Menasco 3-braids were studied in \cite{Pla}. It was shown that 
these double covers are contactomorphic; note, however, that there is a gap in the proof of \cite[Theorem 5]{Pla}. 
\end{remark}

%%%%%%%%%%%%% Ng-Ozsvath-Thurston  %%%%%%%%%%%%%%%%%% 

\subsection{Ng--Ozsv\'ath--Thurston examples}\label{subsec-not} 

In \cite{NgOT}, transverse knots are given as push-offs of Legendrian knots, and the latter 
are represented by grid diagrams of their (smooth) mirrors. We recall how  to obtain a positive transverse push-off 
of a Legendrian knot given by such a grid diagram (cf. \cite{NgOT}). First, let the horizontal segments in the diagram go over the vertical segments (this is opposite to the convention for grid diagrams and produces a front projection 
for the Legendrian knot). Then keep every vertical segment oriented upwards (i.e. has O above X),  
 and replace every vertical segment oriented downwards by the complementary vertical segment. 
The result is a braid that goes from the bottom to the top of the diagram and represents the positive push-off of the given Legendrian knot.  
To obtain the braid for the negative transverse push-off, reverse the 
orientation of the Legendrian knot (by replacing O's by X's and vice versa in the grid diagram), 
and repeat the above procedure.

\begin{figure}[htb!]
\begin{center} 
\includegraphics[scale=1.0]{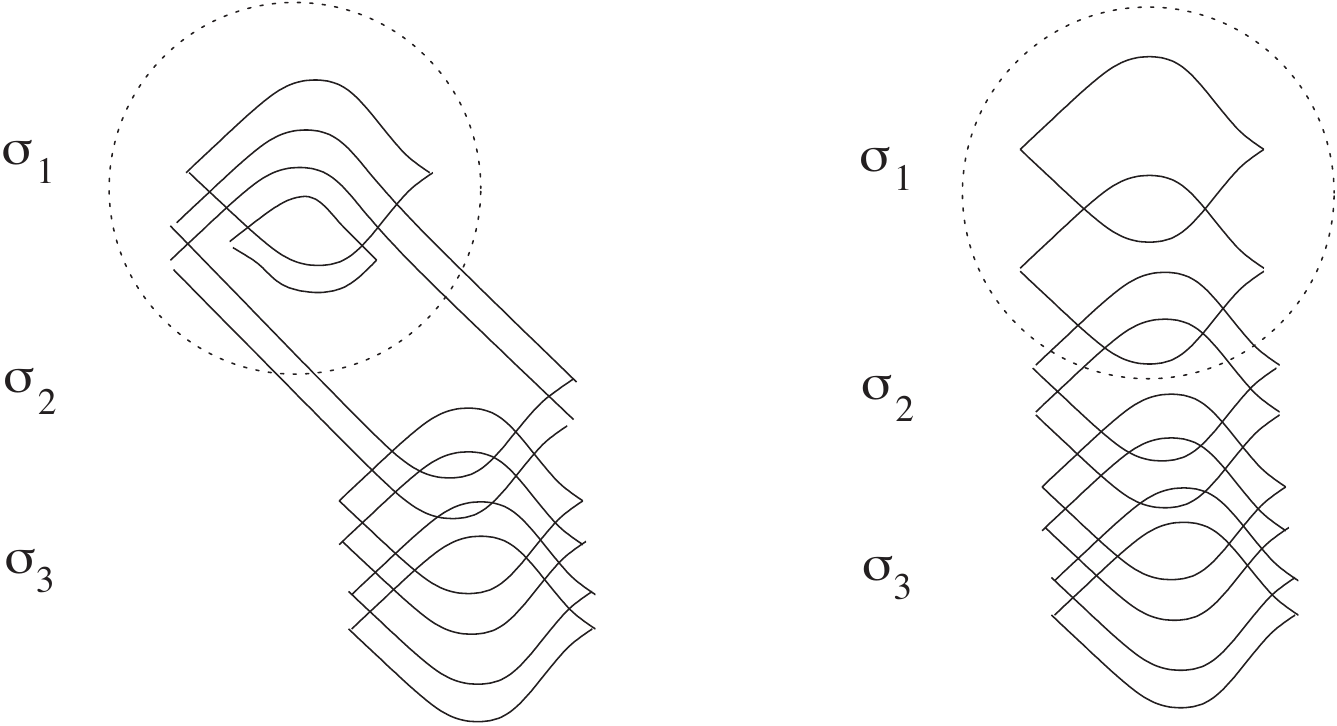}
\caption{Branched double covers of  $L_1^+$ (left) and $L_2^+$ (right).}\label{fig4ngozth}
\end{center}
\end{figure}  

\begin{figure}[htb!]
\begin{center} 
\includegraphics[scale=0.6]{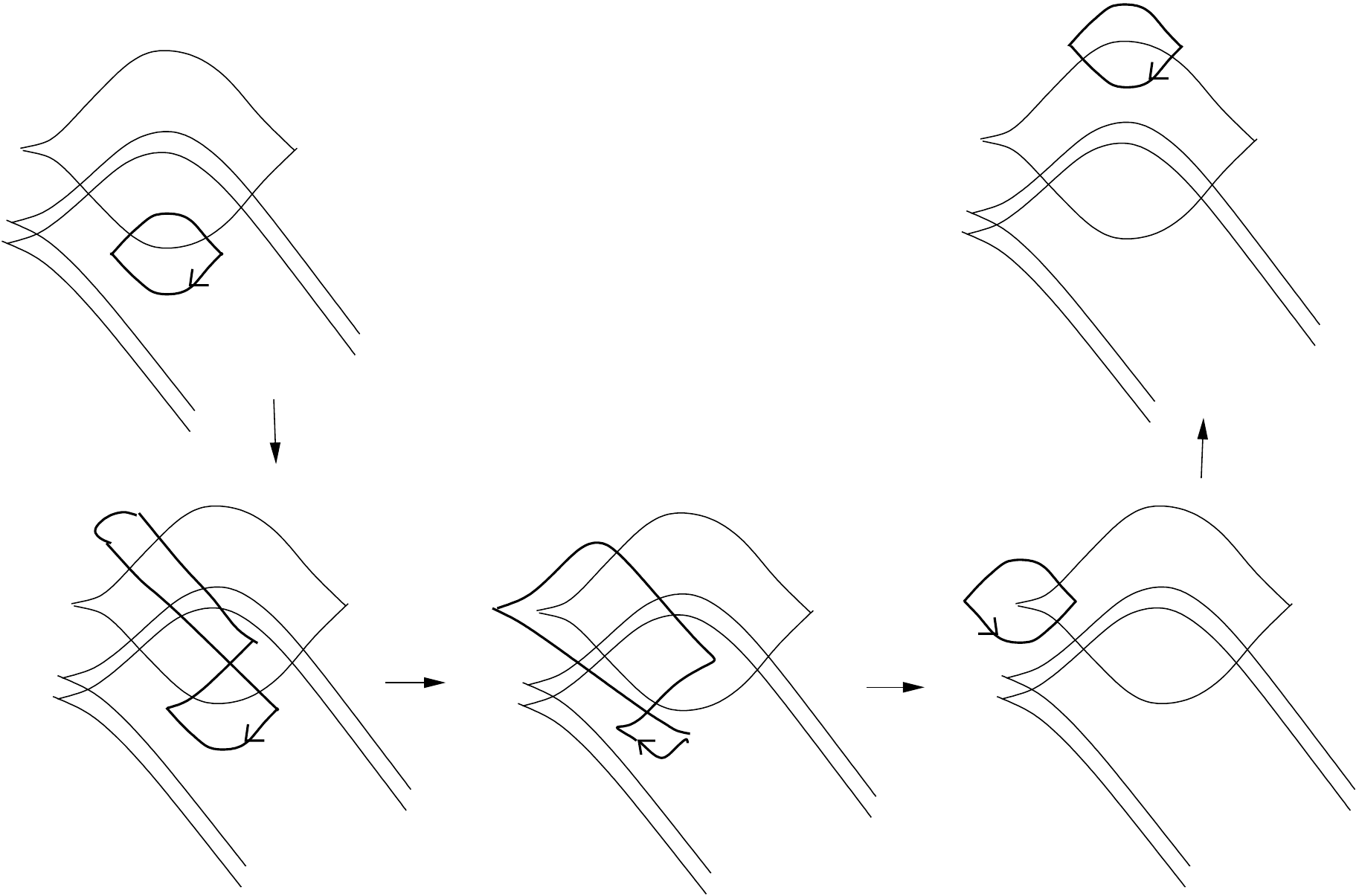}
\caption{The Legendrian links in Figure \ref{fig4ngozth} can be related by Legendrian Reidemeister moves.}   \label{Rmoves}
\end{center}
\end{figure}

We consider transverse push-offs $L_1^+$ and $L_2^+$ of the Legendrian representatives of the pretzel knot $P(-4, -3, 3)$, \cite[Figure 4]{NgOT}. These are given by transverse closed braids 
$$
L_1^+ = \sigma_3^{-1} \sigma_2 \sigma_3 \sigma_1 \sigma_1 \sigma_3 \sigma_2^{-1} \sigma_1 \sigma_2 \sigma_1^{-2}
\quad \text{and} \quad 
L_2^+ = \sigma_3 \sigma_2 \sigma_1 \sigma_3^{-1} \sigma_1 \sigma_2^{-1} \sigma_1 \sigma_2 \sigma_1^{-2} \sigma_3.
$$   
A braid isotopy takes these braids to 
$$
L_1^+ = (\sigma_2 \sigma_1  \sigma_3) \sigma_3 \sigma_1 \sigma_2^{-1} \sigma_3 \sigma_2 \sigma_3^{-2} \sigma_1^{-1}
\quad \text{and} \quad 
L_2^+ = (\sigma_1 \sigma_2 \sigma_3) \sigma_3 \sigma_1^{-1} \sigma_2^{-1} \sigma_3 \sigma_2 \sigma_3^{-2} \sigma_1.
$$
We can now draw surgery diagrams for the double branched covers of $L_1^+$  and $L_2^+$; they are shown on 
Figure~\ref{fig4ngozth}. Recall Remark~\ref{rem-different-choice} and Figure~\ref{AnotherLegReal}. Note that the two surgery links differ only in the circled region; this corresponds to the fact that the braids for $L_1^+$  and $L_2^+$ differ 
only by exchanging two generators $\sigma_1^{-1}$ and   $\sigma_1$ (together with a different choice of 
the open book). We observe that the surgery links are in fact Legendrian isotopic. The isotopy can be performed 
via a sequence of Legendrian Reidemeister moves indicated on Figure~\ref{Rmoves}. 

The transverse push-offs $(L'_1)^+$ and $(L'_2)^+$ of the Legendrian representatives of the 
pretzel knot  $P(-6, -3, 3)$, \cite[Figure 5]{NgOT} can be treated in the same way. Indeed, 
these are given by braids 
\begin{align*}
(L_1')^+&= \sigma_4^{-1} \sigma_3 \sigma_4 \sigma_2 \sigma_1 \sigma_4 \sigma_2 \sigma_1 \sigma_2 \sigma_3^{-1} \sigma_2 \sigma_3 \sigma_2^{-1} \sigma_1^{-1} \sigma_2^{-1} \sigma_1^{-1}, \\  
(L_2')^+&= \sigma_4 \sigma_3 \sigma_2 \sigma_1 \sigma_4^{-1} \sigma_2 \sigma_1 \sigma_2 \sigma_3^{-1} \sigma_2 \sigma_3 \sigma_2^{-1} \sigma_1^{-1} \sigma_2^{-1} \sigma_1^{-1} \sigma_4
\end{align*}
braid isotopic to 
\begin{align*}
(L_1')^+&= (\sigma_2 \sigma_1 \sigma_3 \sigma_4) \sigma_1 \sigma_3 \sigma_4 \sigma_3 \sigma_2^{-1} \sigma_3 \sigma_2 \sigma_3^{-1} \sigma_4^{-1} \sigma_3^{-1} \sigma_4^{-1} \sigma_1^{-1}, \\  
(L_2')^+&= ( \sigma_1 \sigma_2 \sigma_3 \sigma_4) \sigma_1^{-1} \sigma_3 \sigma_4 \sigma_3 \sigma_2^{-1} \sigma_3 \sigma_2 \sigma_3^{-1} \sigma_4^{-1} \sigma_3^{-1} \sigma_4^{-1} \sigma_1.
\end{align*}

As in the previous example, we make a different choice of the initial unknot, and then switch the two factors 
of $\sigma_1$ and $\sigma_1^{-1}$ to relate the braids. The surgery diagrams are very similar to the previous case;
the surgery links have more surgery components, but differ only in the circled region exactly as above, 
and can be related by a sequence of Reidemeister moves. 

It is conjectured in \cite{NgOT} that all pretzel knots $P(-2n, -3, 3)$ are not transversely simple, and if $(L_1^n)^+$, 
$(L_2^n)^-$ 
are the Legendrian representatives  of  $P(-2n, -3, 3)$  similar to those considered above, then  $(L_1^n)^+$ and 
$(L_2^n)^-$ are not transversely isotopic. Our argument, however, clearly generalizes to show that 
the corresponding branched double covers are contactomorphic. 

Moreover, our argument for the knots $L_1^+$ and $L_2^+$ will work for any  
two braids of the form 
$$
K_1= \sigma_1^m \sigma_2 \sigma_1^{-1} w \quad \text{and} \quad K_2= \sigma_1^{-1} \sigma_2 \sigma_1^m w   ,
$$
where $w$ is any braid word on generators $\sigma_2, \dots, \sigma_{n-1}$, and $m>0$. 
Indeed, such two closed braids are isotopic to 
$$
K_1= ( \sigma_1 \sigma_2  \sigma_3 \dots \sigma_{n-1}) \sigma_1^{-1} w' \sigma_1^{m-1} \quad \text{and} \quad
K_2= ( \sigma_2 \sigma_1 \sigma_3  \dots \sigma_{n-1}) \sigma_1^{m-1} w'  \sigma_1^{-1}
$$
and the corresponding surgery diagrams differ by the same local change as above, except that 
instead of the single Legendrian unknot to be moved we have $(m-1)$ copies of its Legendrian push-offs.
A similar sequence of Reidemeister moves can be used to perform this local change. We observe that 
two such braids are in fact related by a negative flype. We thus have 

\begin{prop}\label{K1-K2}
Let $K_1$, $K_2$ be two transverse braids related by the special kind of negative flype satisfying
$$ 
K_1= \sigma_1^m \sigma_2 \sigma_1^{-1} w, \qquad K_2= \sigma_1^{-1} \sigma_2 \sigma_1^m w,
$$
where $w$ is a word in $\sigma_2, \dots, \sigma_{n-1}$, and $m$ is an integer.
Then the branched double covers of $K_1$ and $K_2$ are contactomorphic.
\end{prop} 

\begin{proof} The case $m>0$ is considered above. When $m \leq 0$, the $p$-fold cyclic branched covers for $K_1$ and $K_2$ of any $p$ are overtwisted by Proposition \ref{neg-crossings}. Since they share the homotopy invariants, thus they are contactomorphic.   
\end{proof} 

More generally we have the following.

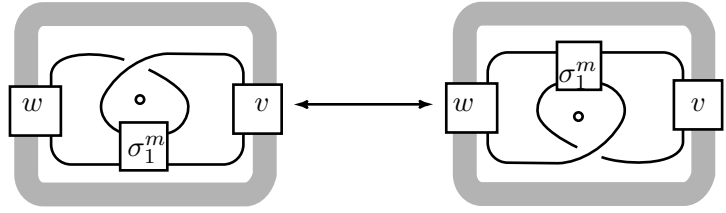
\begin{figure}[htpb] \begin{center}
\begin{picture}(265, 80)
\put(0,0){\includegraphics{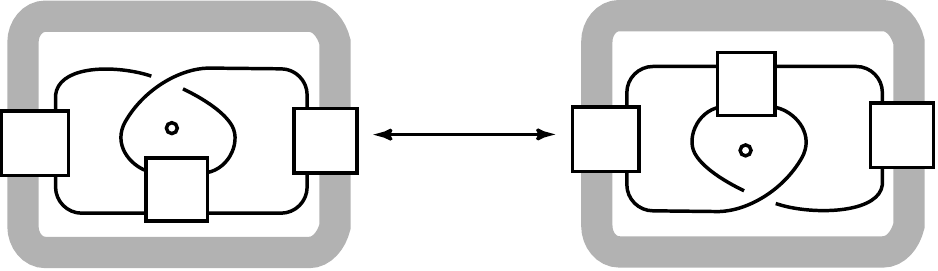}}
\put(93, 37){$v$} \put(258, 37){$v$}
\put(5, 37){$w$} \put(168, 37){$w$}
\put(45, 20){$\sigma_1^m$} \put(208, 48){$\sigma_1^m$}
\end{picture}
\caption{A negative flype move. The gray band means non-braided $(n-2)$ strands.}\label{neg-flype} 
\end{center}\end{figure}

\begin{theorem}\label{neg-flype-thm}
Suppose $K_1$ and $K_2$ are related each other by a negative flype move sketched in Figure~\ref{neg-flype}, i.e., 
$$K_1 = \sigma_1^m v \sigma_1^{-1} w, \qquad K_2= \sigma_1^{-1} v \sigma_1^m w,$$
where $v$ and $w$ are any braid words in generators $\sigma_2, \dots, \sigma_{n-1}$ and $m \in \ZZ$. Then the branched double covers $(\Sigma_2(K_1), \xi_2(K_1))$ and $(\Sigma_2(K_2), \xi_2(K_2))$ are contactomorphic.
\end{theorem}

\begin{proof}
Consider positive stabilizations of $K_1$ and $K_2$. Since a positive stabilization preserves transverse knot type, we use the same notations $K_1, K_2$. Let $v'$ (resp. $w'$) be the braid words in $\sigma_3, \dots, \sigma_n$ obtained from $v$ (resp. $w$) by translation $\sigma_k \mapsto \sigma_{k+1}$. Then we have
\begin{eqnarray*}
K_1 &=&  \sigma_1^m \sigma_1^{-1} \sigma_1 v \sigma_1^{-1} w  \qquad \mbox{isotopy}\\
       &=& \sigma_2^m  \sigma_2^{-1} \underline{\sigma_1} \sigma_2 v' \sigma_2^{-1} w' \qquad (+)\mbox{stabilization}\\
       &=&  \sigma_2 \sigma_1^m v'  \sigma_2  \sigma_1^{-1}  \sigma_2^{-1} w'   \qquad \mbox{isotopy}\\
       &=& \sigma_2 v' \underline{\sigma_1^m \sigma_2  \sigma_1^{-1}}  \sigma_2^{-1} w'  \qquad \mbox{isotopy.}
\end{eqnarray*}       
Similarly we have 
$$K_2 = \sigma_2  \sigma_1^{-1}  v' \sigma_2  \sigma_1^m \sigma_2^{-1} w' = \sigma_2  v' \underline{\sigma_1^{-1} \sigma_2  \sigma_1^m} \sigma_2^{-1} w'.$$
Thus they satisfy the condition of Proposition~\ref{K1-K2}.
\end{proof}

\begin{example}
Let $L_1, L_2$ (resp. $L_1', L_2'$) be the Legendrian $m(10_{132})$ (resp. $m(12n_{200})$) knots studied in \cite{NgOT}. Let $M_1, M_2$ be the Legendrian $(2, 3)$-cables of the $(2, 3)$-torus knot  found in \cite{EH1} \cite{MM}. The positive push-offs of every pair satisfy the condition of Theorem~\ref{neg-flype-thm}. Therefore, double branched covers for each pair are contactomorphic.
\end{example}

\begin{proof}
It is shown in \cite{NgOT} that the closed braids $(L_1)^+, (L_2)^+$ (resp. $(L_1')^+, (L_2')^+$) only differ in the dashed boxes sketched in Figure~\ref{neg-flype-1} and are related to each other by a negative flype.
\begin{figure}[htpb] \begin{center}
\begin{picture}(333, 206)
\put(0,0){\includegraphics{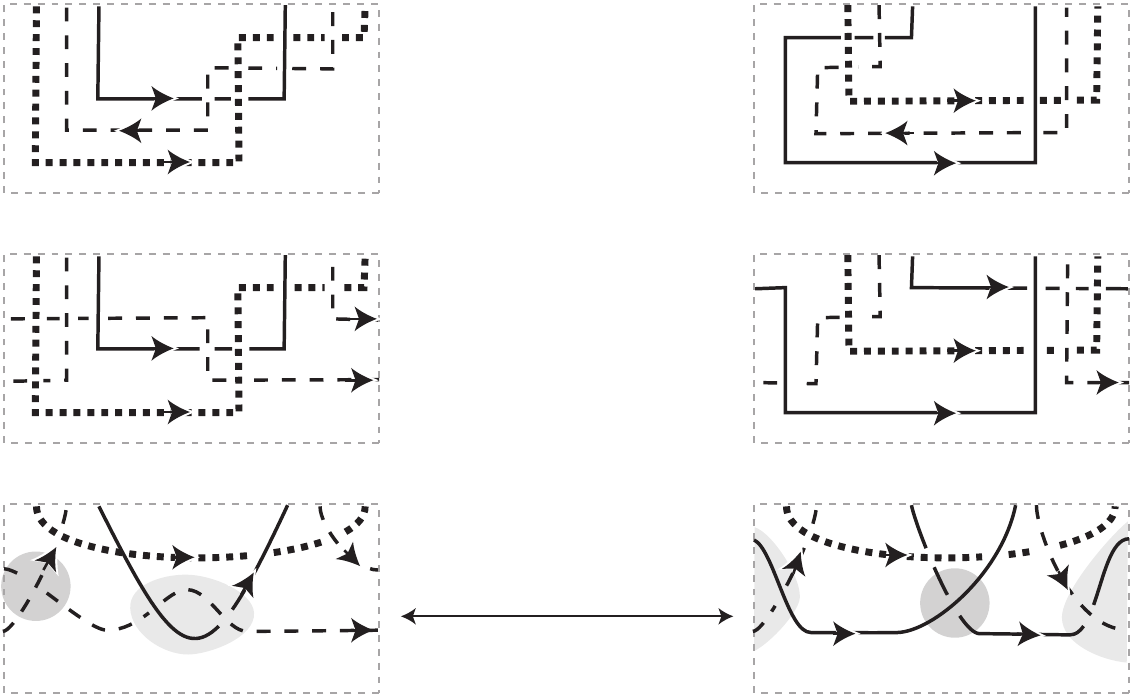}}
\put(130, 30){Negative flype}
\put(95, 150){$L_1$} \put(310, 150){$L_2$}
\end{picture}
\caption{}\label{neg-flype-1} 
\end{center}\end{figure}  

Similarly, the closed braids $(M_1)^+, (M_2)^+$ only differ in the dashed boxes sketched in Figure~\ref{neg-flype-2} and are related to each other by a negative flype. (Here we use the Legendrian fronts for $M_1, M_2$ given in \cite{NgOT}.) 
\begin{figure}[htpb] \begin{center}
\begin{picture}(360, 170)
\put(0,0){\includegraphics{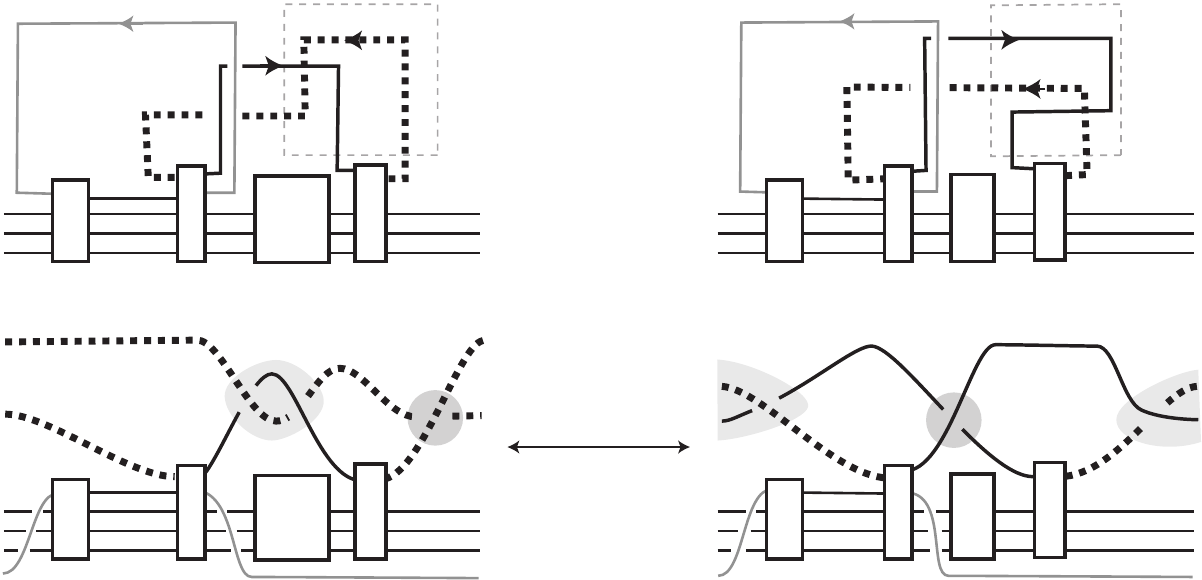}}
\put(155, 28){Negative}
\put(165, 18){flype}
\put(20, 130){$M_1$} \put(223, 130){$M_2$}
\end{picture}
\caption{}\label{neg-flype-2} 
\end{center}\end{figure}
\end{proof}

%%%%%%%%%%%%%%%%%%%%%%%%%%%%%%%%%

\subsection{Heegaard Floer contact invariants} \label{OSi}

We are unable to determine whether the branched double covers for some of the examples in  \cite{NgOT} are contactomorphic. However, we can show that the Heegaard Floer contact invariants \cite{ContOS} also fail to distinguish between these contact manifolds.

\begin{prop} Let $K_1, K_2$ be one of the pairs $L_1^+, L_1^-$ or $(L'_2)^+, (L'_2)^-$ of transverse push-offs of Legendrian pretzel knots $P(-4, -3, 3)$ and $P(-6, -3, -3)$ 
\cite[Figure 4, 5]{NgOT}.  
Then $c(\xi_2(K_1))= c(\xi_2(K_1))$.
\end{prop}

\begin{proof} We use the following result of Lawrence Roberts \cite{Ro}, conjectured in \cite{Pla}.
For a transverse link $L$, let $\psi(L)$ be its transverse  invariant in reduced Khovanov homology, 
\cite{PlKh}. Recall that by \cite{OSdouble}, there is a spectral sequence that relates the reduced 
Khovanov homology of $L$ to the Heegaard Floer homology $\h{HF}(-\Sigma_2(L))$ (when one works with 
coefficients in $\ZZ/2$. Then the element 
 $\psi(L)$  canonically corresponds to the Heegaard Floer contact invariant under this spectral sequence $c(\xi_2(L))$. 

Our result now follows from the fact that $\psi(K_1)=\psi(K_2)$ for all examples in question \cite{NgOT}.
Indeed,  for these pretzel knots Khovanov 
homology in the relevant bi-degree has rank 1, and  $\psi(K_1)=\psi(K_2)\neq 0$ because all braids are quasipositive \cite{PlKh}.  
\end{proof}

\end{document}